\documentclass[11pt, a4paper]{amsart}

\usepackage[colorlinks, linkcolor=blue, anchorcolor=blue, citecolor=green]{hyperref}

\usepackage{graphics,epic}
\usepackage{amsmath,amssymb, amsthm,mathrsfs}
\usepackage[all,2cell]{xy}
\usepackage{shorttoc}
\usepackage{tikz}
\usepackage{tikz-cd}

\newtheorem{theorem}{Theorem}[section]
\newtheorem*{theorem*}{Theorem}

\newtheorem*{conjecture*}{Conjecture}

\newtheorem{definition}[theorem]{Definition}

\newtheorem{thm}[theorem]{Theorem}
\newtheorem{lem}[theorem]{Lemma}
\newtheorem{prop}[theorem]{Proposition}
\newtheorem{cor}[theorem]{Corollary}

\newtheorem{rk}[theorem]{Remark}

\newcommand{\ie}{{\em i.e.}\ }

\newcommand{\opname}[1]{\operatorname{\mathsf{#1}}}

\renewcommand{\mod}{\opname{mod}\nolimits}

\newcommand{\proj}{\opname{proj}\nolimits}
\newcommand{\inj}{\opname{inj}\nolimits}
\newcommand{\dimv}{\underline{\dim}\,}

\newcommand{\add}{\opname{add}\nolimits}

\newcommand{\ind}{\opname{ind}}

\renewcommand{\P}{\mathbb{P}}

\newcommand{\ra}{\rightarrow}

\newcommand{\sgn}{\opname{sgn}}

\newcommand{\Hom}{\opname{Hom}}

\newcommand{\Ext}{\opname{Ext}}

\newcommand{\End}{\opname{End}}

\newcommand{\ten}{\otimes}
\newcommand{\lten}{\overset{\mathbf{L}}{\ten}}

\newcommand{\md}{\mathcal{D}}
\newcommand{\mf}{\mathcal{F}}

\newcommand{\mt}{\mathcal{T}}

\newcommand{\mx}{\mathcal{X}}
\newcommand{\my}{\mathcal{Y}}

\newcommand{\tauni}{\tau^{-1}}

\newcommand{\xra}{\xrightarrow}

\begin{document}

\title{Functor-induced isomorphisms and $G$-matrices}\thanks{Partially supported by the National Natural Science Foundation of China (Grant No.   12471037) }

\author[Shengfei Geng]{Shengfei Geng}
\address{Shengfei Geng\\Department of Mathematics\\SiChuan University\\610064 Chengdu\\P.R.China}
\email{genshengfei@scu.edu.cn}


\begin{abstract}
In this paper, we explore how functor-induced isomorphisms are encoded by $G$-matrices. We first show that the Grothendieck group isomorphism induced by a tilting module can be realized via  the $G$-matrix  of this tilting module. Building on this, we compare $g$-vectors for a tilted algebra and its associated hereditary algebra, and provide $G$-matrix interpretations of the Coxeter transformation, the Nakayama functor, and the Auslander-Reiten translation for suitable algebras. Furthermore, we demonstrate the geometric representations of any Weyl group can be realized by the $G$-matrices of tilting modules or support $\tau$-tilting modules arising from generalized preprojective algebras. Consequently, the group generated by the transpose of these $G$-matrices is isomorphic to the corresponding Weyl group.  Finally, we show that the Grothendieck group isomorphism induced by a $2$-term silting complex can also be realized via the $G$-matrix of this $2$-term silting complex.

\end{abstract}
\maketitle
\tableofcontents

\section{Introduction}

In the theory of cluster algebras (cf. \cite{FZ02,FZ03,FZ07}), \( g \)-vectors serve as fundamental invariants that encode the grading of the cluster variables. An important result, established by \cite{GHKK18}, asserts that distinct cluster variables—and more generally, different cluster monomials—possess distinct \( g \)-vectors. For a given cluster, the associated \( G \)-matrix is formed by taking the \( g \)-vectors of its cluster variables as columns. This matrix exhibits several key properties: it is unimodular (with determinant \( \pm 1 \)), sign-coherent (each row is either entirely non-negative or entirely non-positive), and dual to the \( C \)-matrix (cf. \cite{GHKK18}). In particular, for cluster algebras of finite type, the $g$-vectors of all cluster variables constitute the vertex set of a polytope known as the generalized associahedron. The normal fan of this polytope forms a complete simplicial fan in $\mathbb{R}^n$, reflecting the rich combinatorial structure of the cluster algebra.

In the categorification of cluster algebras, \( g \)-vectors are similarly defined for modules over finite dimensional algebras, where they play an important role in the study of \( \tau \)-tilting theory. Notably, distinct \( \tau \)-rigid pairs are characterized by distinct \( g \)-vectors. The \( G \)-matrix of a basic \( \tau \)-tilting pair—whose columns are the \( g \)-vectors of its indecomposable direct summands—forms a \( \mathbb{Z} \)-basis of \( \mathbb{Z}^n \), and is unimodular and row sign-coherent. This parallel underscores the deep interplay between cluster algebras and their categorical counterparts. For further details and results, we refer to  \cite{DK08, AIR14, DIJ19,BST19,A21,C21,PYK23, AHIKM} and so on.

It is known that  certain functors  can induce isomorphisms between Grothendieck groups of relevant  categories. A natural question is whether these functor-induced isomorphisms can be described by \( G \)-matrices. In this paper, we explore how these functor-induced isomorphisms are encoded by $G$-matrices.

Let $k$ be an algebraically closed field,  $A$ be a finite dimensional $k$-algebra.  For an \( A \)-module \( M \), denote its \( g \)-vector by \( g^M \). For a basic tilting module \( T = \bigoplus_{i=1}^n T_i \), its $G$-matrix \( G_T \) is the matrix \( (g^{T_1}, \cdots, g^{T_n}) \). Let $B=\End_A T$, it is known that  the correspondence 
      $$g:\dimv N\to \dimv\Hom_A(T,N)-\dimv\Ext_A^1(T,N),$$
     where $N$ is a $A$-module, induces an isomorphism \( g: K_0(A) \to K_0(B) \) between the Grothendieck groups of \( A \) and \( B \). Our first observation is
\begin{thm}(Theorem~\ref{t: iso induced by G_T})
For any $A$-module $N$, we have $$G_T^t\dimv N=\dimv\Hom_A(T,N)-\dimv\Ext_A^1(T,N).$$
where   \( G_T^t \) is the  transpose of  \( G_T \).
\end{thm}
Building on this result, we present several applications. First, we reformulate \cite[Corollary 7.6]{GLS06} in terms of the $G$-matrix and give a $G$-matrix interpretations of  the mutation on certain kinds of matrices (see Theorem~\ref{t:mutation via g-matrix}). Then we compare the $g$-vectors of a tilted algebra with its associated hereditary algebra (see Proposition~\ref{p: g-vectors between hereditary algebras and tilted algebras}). Moreover, recalling that $DA$ is a tilting $A$-module if and only if $A$ is $1$-Gorenstein (where $D=\Hom_k(-,k)$), we show that  the Coxeter matrix $\Phi_A$, the Nakayama functor $\nu$ and the Auslander–Reiten translation $\tau$ can be described by $G_{DA}=(g^{I(1)},\cdots,g^{I(n)})$ in certain cases, where $I(1),\cdots, I(n)$ are all the indecomposable injective $A$-modules.
 \begin{prop}
     \begin{itemize}
         \item [(1)](Proposition~\ref{p:judge 1-gorenstien alg}) If $A$ is a $1$-Gorenstein algebra and $\mathbf{det} \mathbf{C}_A\neq 0$,  then $\Phi_A=-(G_{DA}^{-1})^t$, where $\mathbf{C}_A$ the Cartan matrix of $A$.
         \item [(2)](Proposition~\ref{p: nakayama for self-injective algebra}) If $A$ is a self-injective algebra, then for any $A$-module $X$, we have $\dimv\nu(X)=G_{DA}\dimv X=(G_{DA}^{-1})^t\dimv X$ and $g^{\nu X}=G_{DA}g^X=(G_{DA}^{-1})^tg^X$.
         \item[(3)](Proposition~\ref{p: tau induced by g-vector for hereditary alg}) If $A$ is a hereditary algebra, then for any indecomposable non-projective $A$-module $X$, we have $\dimv \tau X=-(G_{DA}^{-1})^t\dimv X$ and $g^{\tau X}=-G_{DA}g^X$.
     \end{itemize}
 \end{prop}
Furthermore, when \( A = kQ/I \) is the Auslander algebra of \( k[x]/(x^n) \),  following from \cite{IZ20}, there exists a bijection between basic tilting \( A \)-modules and the symmetric group \( \mathfrak{S}_n \). Thus, for each word \( w \in \mathfrak{S}_n \), there corresponds a tilting \( A \)-module \( I_w \), which is simultaneously a tilting ideal of \( A \).
Let \( e_1, \cdots, e_n \) be the idempotents associated to the vertices of \( Q \), and fix  
\[
G_{I_w} = \left( g^{e_1 I_w}, \cdots, g^{e_n I_w} \right).
\]
Note that \( w \) can also be viewed as an element of \( \mathfrak{S}_{n+1} \). Denote by \( R_w \) the geometric representation matrix of \( w \) in \( \mathfrak{S}_{n+1} \) (see Subsection~\ref{s:def of symmetric group}). Then based on the work of \cite{IZ20}, we have
\begin{thm}(Theorem~\ref{t: symmetric group via g-matrices})
  Let  $A$ be the Auslander algebra of $k[x]/ (x^n)$, then for  each word $w\in \mathfrak{S}_n$, we have $G_{I_w}^t=R_w$.
  \end{thm}
This result can be generalized to any Weyl group. Let $C$ be a symmetrizable generalized Cartan matrix with a symmetrizer $D$. According to  \cite{GLS17}, there is an associated generalized preprojective algebra $\Pi=\Pi(C,D)$. When $C$ is symmetric and $D$ is the identity matrix, $\Pi$ reduces to the classical preprojective algebras for acyclic quivers.
Let $\overline{Q}
$ be the  quiver of $\Pi$ with
vertex set $Q_0 := \{1, \cdots , n\}$.  For each $i\in Q_0,$ denote by $I_i$ the two-sided ideal $\Pi(1-e_i)\Pi$ and
 $\langle I_1,I_2,\cdots,I_n\rangle$  the ideal semigroup generated by $I_1,I_2,\cdots,I_n.$  Denote by $W(C)$  the corresponding Weyl group.
For a reduced expression $w=s_{i_1}\cdots s_{i_l}\in W(C)$, denoted by $I_w=I_{i_1}\cdots I_{i_l}$. If $C$ has no component of Dynkin type, $I_w$ is a cofinite tilting ideal of $\Pi$ and $\End_{\Pi}I_w\cong \Pi$ (cf. \cite{BIRS09, FG19}). If  $C$ is of Dynkin type, $I_w$ is a support $\tau$-tilting module over $\Pi$ (cf. \cite{BIRS09, FG19}). Fix $G_{I_w}$ as Subsection~\ref{s:reflection via g-matrix},
 denote by $R_w$ the geometric representation matrix of $w$, and by $\Sigma_w$ the contragredient matrix of $w$ (see Subsection~\ref{ss:Weyl group}).
 When $\Pi$ is a classical preprojective algebra of extended Dynkin type, it is shown that $G_{I_w}=\Sigma_{w^{-1}}$ (cf. \cite{IR08}), while if $\Pi$ is a classical preprojective algebra of  Dynkin type, it is shown that $G_{I_w}=\Sigma_{w^{-1}}$ (cf. \cite{M14}).
Motivated by these results, we get
 \begin{thm}(Theorem~\ref{t: Weyl group via g-matrices})
      Let $C$ be a symmetrizable Cartan matrix with a symmetrizer $D$ and $\Pi=\Pi(C,D)$ the generalized preprojective algebra. Then for any  word $w\in W(C)$, we have $$G_{I_w}=\Sigma_{w^{-1}}=R_w^t.$$
     In particular, the group $\mathcal{G}_{\Pi}$ generated by $\{G_{I_w}^t \mid w\in W(C)\}$ is isomorphic to the Weyl group $W(C)$.
 \end{thm}

Finally, let $A$ be a finite dimensional  algebra. 
We consider 2-term silting complexes \( P^* \) in \( K^b(\mathsf{proj}\, A) \).
Since there is a bijection between 2-term silting complexes and $\tau$-tilting pair (cf. \cite{AIR14}), so we also have a $G$-matrix $G_{P^*}$ for any 2-term silting complex $P^*$.
Let $B=\End_{K^b(\proj A)} P^*$. By \cite{Hu22},  the correspondence 
      $$f:\dimv X\to \dimv \Hom_{\md^b(A)}(P^*,X)-\dimv \Hom_{\md^b(A)}(P^*,X[1]),$$
      where $X$ is an $A$-module, induces an isomorphism $f:K_0(A)\to K_0(B)$ of the Grothendieck groups of $A$ and $B$.  
Then we can give a generalization of Theorem~\ref{t: iso induced by G_T}.
\begin{thm}(Theorem~\ref{t:iso induced by a 2-term silting complex})
For each $A$-module $X$, we have
$$G_{P^*}^t\dimv X= \dimv \Hom_{\md^b(A)}(P^*,X)-\dimv \Hom_{\md^b(A)}(P^*,X[1]).$$

\end{thm}

\subsection*{Structure} The structure of the paper is outlined as follows. In Section~\ref{s:g-vectrs and G-matrix}, we recall some basic definitions and results of $g$-vectors and $G$-matrices. In Section~\ref{s:Tilting-induced Grothendieck group isomorphism  and G-matrix}, we show that the Grothendieck group isomorphism induced by a tilting module can be realized by the transpose of the $G$-matrix of this tilting module. Then we provide some applications,  we reformulate \cite[Corollary 7.6]{GLS06} in terms of the $G$-matrix, and compare the $g$-vectors of a tilted algebra with its associated hereditary algebra. In Section~\ref{s:other functor-induced iso and G-matrix}, we provide $G$-matrix interpretations of the Coxeter transformation, the Nakayama functor, and the Auslander–Reiten translation for suitable algebras.
In Section~\ref{s:Symmetric groups via G-matrices} and Section~\ref{s:Weyl groups via G-matrices}, we demonstrate that each geometric representation matrix of any element from the symmetric group or from the Weyl group can be expressed as the transpose of the $G$-matrix of some tilting module or support $\tau$-tilting module. In Section~\ref{s: iso induced by 2-term silting complex and G-matrix}, we show that the Grothendieck group isomorphism induced by a $2$-term silting complex can also be realized via the transpose of the $G$-matrix of this $2$-term silting complex.

\subsection*{Convention}

Throughout this paper, let $k$ be an algebraically closed field. By a finite dimensional algebra, we always mean a basic finite dimensional algebra over $k$.  
For such an algebra \(A\), 
\begin{itemize}
    \item[-] \(\bmod A\) is the category of finitely generated right \(A\)-modules;
\item[-] \(\ind A\) is the set of isomorphism classes of indecomposable  \(A\)-modules;
    \item[-] \(\tau\) is its Auslander–Reiten translation; 
    \item[-] \(\{e_{1},\dots ,e_{n}\}\) is a fixed complete set of pairwise orthogonal primitive idempotents;  
    \item[-] \(P(i)=e_{i}A\),  \(S_{i}=\operatorname{top}P(i)\),  \(I(i)=D(Ae_{i})\) are the corresponding indecomposable projective, simple and injective  modules, with \(D=\operatorname{Hom}_{k}(-,k)\);
\end{itemize}
For \(M\in\bmod A\),
\begin{itemize}
    \item[-] \(|M|\) counts the non-isomorphic indecomposable summands of \(M\); 
    \item[-] \(\operatorname{proj}_{\!A}M\) and \(\operatorname{inj}_{\!A}M\) are its projective and injective dimensions;  
    \item[-] \(\operatorname{add} M\) is the subcategory of direct summands of finite direct sums of  \(M\).
\end{itemize}
We write
\begin{itemize}
    \item[-] \(\md^{b}(\bmod A)\) for the bounded derived category of $\mod A$;  
    \item[-] \(K_{0}(A)\) for the Grothendieck group of \(\bmod A\) (or equivalently of \(\md^{b}(\bmod A)\));  
    \item[-] \([M]\) for the class of \(M\in\md^{b}(\bmod A)\) in \(K_{0}(A)\).
   \end{itemize}
For any matrix \(G\), \(G^{t}\) is its transpose; if \(G\) is square, \(\mathbf{det} G\) is its determinant, and \(\mathbf{I}\) denotes the identity matrix of appropriate size.

\section{On g-vectors and G-matrices}\label{s:g-vectrs and G-matrix}
\subsection{Tilting modules and $\tau$-tilting pairs}
Let \( A \) be a finite dimensional \( k \)-algebra.  
\begin{definition}
  A module \( T \in \mathsf{mod}\,A \)  is called  a tilting module if it satisfies the following three conditions:
    \begin{itemize}
        \item[(T1)]$\proj_A T\leq 1;$
        \item[(T2)]$\Ext_A^1(T,T)=0;$
        \item[(T3)]There exists a short exact sequence $0\ra A\ra T_0\ra T_1\ra 0$ for some $T_0,T_1\in\add T$.
    \end{itemize}
\end{definition}
A module satisfying (T1) and (T2) is said to be {\it partial tilting}. It is known that a partial tilting module \( T \) is a tilting module if and only if \( |T| = |A| \). 

\begin{definition}Let $M$ be a $A$-module and $P$ be a  projective $A$-module.
\begin{itemize}
    \item [(1)]$M $ is  a {\it $\tau$-rigid} $A$-module if $\Hom_{A}(M,\tau M) = 0.$
    \item[(2)] $M$ is a {\it $\tau$-tilting} $A$-module if it is $\tau$-rigid and $|M| = |A|.$
    \item[(3)]  $M$ is a {\it support $\tau$-tilting} $A$-module if there exists an idempotent $e \in A$ such
 $M$ is a $\tau$-tilting $(A/\langle e\rangle)$-module.
  \item[(4)] A pair $(M,P)$ is a {\it $\tau$-rigid pair} if $M$ is $\tau$-rigid and
$\Hom_{A}(P,M) = 0.$
    \item[(5)] A $\tau$-rigid pair $(M,P)$ is a {\it $\tau$-tilting}  pair if $|M| + |P| = |A|$.
\end{itemize}
\end{definition} 
 For a basic $\tau$-tilting pair $(M,P)$, by \cite{AIR14} , $P$ is uniquely determined by $M$. It is clear that a $\tau$-rigid $A$-module $M$ is partial tilting if $\proj_A M\leq 1$. In particular,
 a $\tau$-tilting module $T$ is  tilting  if and only if $\proj_A T\leq 1$.

\subsection{$g$-vector}\label{ss:g-vector}
For an $A$-module $M$, let $$\bigoplus_{i=1}^nP(i)^{b_i}\xrightarrow{f_M}\bigoplus_{i=1}^nP(i)^{a_i}\to M\to 0$$ be the minimal projective presentation of $M$, where $a_i,b_j$ are non-negative integers. The {\it $g$-vector} $g^M_A$ (simply denoted by $g^M$) associated with $M$ is defined as
 \[
 g_A^M=[a_1-b_1,\cdots, a_n-b_n]^t.
 \]
where $t$ means transpose.

For a $\tau$-rigid pair $(M,P)$, {$g^{(M,P)}:=g^M-g^P$} is called the {\it $g$-vector of $(M,P)$}. Let  $(M,P)$ be a basic  $\tau$-tilting pair, written as $M=\bigoplus_{i=1}^t M_i,$ $ P=\bigoplus_{i={t+1}}^n P_i$ with each $M_i$ and $P_j$ indecomposable.  
The \( G \)-matrix of \( (M,P) \) is the integer matrix  {\it $$G_{(M,P)}=(g^{M_1}, \cdots, g^{M_t}, -g^{P_{t+1},}\cdots,-g^{P_n}). $$}
By \cite{AIR14}, the columns of \( G_{(M,P)} \) form a \( \mathbb{Z} \)-basis of \( \mathbb{Z}^{n} \). Hence,  
\[
\mathbf{det} G_{(M,P)}\in\{1,-1\}.
\]  
Consequently, one can define an integer matrix $C_{(M,P)}= (G_{(M,P)}^t)^{-1}$, called the $C$-matrix of $(M,P)$.
Similarly,  the  {$D$-matrix of $(M,P)$} is defined as {\it $$D_{(M,P)}=(\dimv {M_1}, \cdots, \dimv {M_t}, -\dimv {P_{t+1},}\cdots,-\dimv {P_n}).$$}

\begin{thm}\label{t: basic results of g-vector and g-matrix}
Let $A$ be a finite dimensional algebra. Then
    \begin{itemize}
    \item[(1)]  Different $\tau$-rigid $A$-pairs have different $g$-vectors \cite{AIR14}.
    \item[(2)]  Every \( G \)-matrix of $A$ is  row sign-coherent \cite{AIR14} .
    \item[(3)]   Every \( C \)-matrix of $A$ is  column sign-coherent \cite{Fu17}.
    
\end{itemize}
\end{thm}

\subsection{Dimension vectors and $g$-vectors}
Let $A$ be a finite dimensional $k$-algebra. 
Denote by $$\mathbf{C}_A=(\dimv P(1),\cdots,\dimv P(n))$$  the Cartan matrix of $A$.  It is known that $\mathbf{C}_A^t=(\dimv I(1),\cdots,\dimv I(n))$.
For a $A$-module  $M$, take a minimal projective presentation $ P_1\xra{f_1}P_0\ra M\ra 0$ of $M$,
then it is easy to get that 
\begin{eqnarray}\label{g: C_Ag^M=p0-p1}
    \mathbf{C}_Ag^M =\dimv P_0-\dimv P_1.
\end{eqnarray}
Denote by $\nu=D\Hom(-,A)$ the Nakayama functor. Because $\nu P(i)=I(i)$, then one can  get that 
\begin{eqnarray}\label{g: C_A^tg^M=vp0-vp1}
   \mathbf{C}_A^tg^M =\dimv \nu P_0-\dimv \nu P_1.
\end{eqnarray}
Moreover, we have
\begin{lem}\label{l: proj dim and g-vector}
    Let $A$ be a finite dimensional algebra.
    \begin{itemize}
        \item[(1)] For any $A$-module $M$,  $\proj_A M\leq 1$ if and only if $\dimv M=\mathbf{C}_Ag^M $.
        \item[(2)] For any $A$-module $M$,  suppose no direct summand of $M$ is projective, then $\inj_A \tau M\leq 1$ if and only if $\dimv \tau M=-\mathbf{C}_A^tg^M $.
\item[(3)] Let $T$ be a $\tau$-tilting $A$-module, then $T$ is a tilting $A$-module if and only if $D_T=\mathbf{C}_AG_T $.

    \end{itemize}
\end{lem}
\begin{proof}
For a $A$-module  $M$, take a minimal projective presentation $$ P_1\xra{f_1}P_0\ra M\ra 0$$ of $M$.

For $(1)$, if $\proj_A M\leq 1$, it is easy to get that $\dimv M=\mathbf{C}_Ag^M $. On the other hand, if $\dimv M=\mathbf{C}_Ag^M $, by~(\ref{g: C_Ag^M=p0-p1}),  $$\dimv P_0-\dimv P_1=\mathbf{C}_Ag^M =\dimv M,$$
which means $\dimv \mathbf{ker}(f_1)=0$. So $\proj_AM\leq 1$.

For $(2)$, if no direct summand of $M$ is projective, then we  have a minimal injective presentation 
$$0\ra \tau M\ra \nu P_1\xra{\nu (f_1)} \nu P_0$$ of $\tau M$. Thus if  $\inj_A \tau M\leq 1$, it is clear that $\dimv \tau M=-\dimv \nu P_0+\dimv \nu P_1=-\mathbf{C}_A^tg^M $. On the other hand, if $\dimv \tau M=-\mathbf{C}_A^tg^M$, by~(\ref{g: C_A^tg^M=vp0-vp1}), we have $\dimv \tau M=-\mathbf{C}_A^tg^M=-\dimv \nu P_0+\dimv \nu P_1$, which means $\mathbf{coker} (\nu(f_1))=0$, thus $\inj_A \tau M\leq 1$.

 $(3)$ can be get from $(1)$ directly.
    \end{proof}

\begin{cor}\label{c: modules determind by their dimension vectors}
 Let $A$ be a finite dimensional algebra. The following are equivalent:
 \begin{itemize}
\item[(1)] $\mathbf{det} \mathbf{C}_A\neq 0$;
 \item[(2)] different partial tilting $A$-modules have different dimension vectors; 
\item[(3)] for any tilting module $T$, $\mathbf{det} D_T\neq 0$;
      \item [(4)] there is a tilting module $T$ such that  $\mathbf{det} D_T\neq 0$.
     
 \end{itemize}
   \end{cor}
\begin{proof}
Let $T$ be a tilting $A$-module, so $\mathbf{det} G_T\neq 0$. By $\mathbf{C}_AG_T=D_T$, so $\mathbf{det} D_T\neq 0$ if and only if $\mathbf{det} \mathbf{C}_A\neq 0$.
So (1),(3) and (4) are equivalent. 

(2)$\Rightarrow$(1) is clear since projective modules are partial tilting.

(1)$\Rightarrow$(2)
Let $X,Y$ be two partial tilting $A$-modules with  $\dimv X=\dimv Y$. So $X,Y$ are also $\tau$-rigid $A$-modules. Because $\proj_A X\leq 1$ and $\proj_A Y\leq 1$, one can get that  $\dimv X=\mathbf{C}_Ag^X$, $\dimv Y=\mathbf{C}_Ag^Y$. Since $\mathbf{det} \mathbf{C}_A\neq 0$, we have $g^X=g^Y$. By~Theorem~\ref{t: basic results of g-vector and g-matrix},
 different $\tau$-rigid $A$-module have different $g$-vectors, thus we must have $X=Y$.
  \end{proof}

\subsection{Key Lemma} Denote by $\langle-,-\rangle$ the canonical inner product of $\mathbb{R}^n$. The following lemma is useful.
\begin{lem}\cite[Theorem 1.4]{AR85}\label{l:key lemma}
    Let $M,N$ be two $A$-modules, then $$\langle g^M,\dimv N\rangle=\dim\Hom_A(M,N)-\dim\Hom_A(N,\tau M).$$
\end{lem}

\section{Tilting-induced Grothendieck group isomorphisms  and G-matrices}\label{s:Tilting-induced Grothendieck group isomorphism  and G-matrix}
\subsection{Grothendieck group isomorphism induced by a tilting module and the G-matrix}
Let $A$ be a finite dimensional algebra. 
We identify the Grothendieck group \( K_0(A) \) with \( \mathbb Z^n \) via the dimension vector map
\[
\dimv:\mod A\longrightarrow \mathbb Z^n,\qquad S_i\longmapsto \mathbf e_i,
\]
where \( \mathbf e_1,\cdots,\mathbf e_n \) denotes the standard basis of \( \mathbb Z^n \).

Let $T=\bigoplus_{i=1}^n T_i$ be a tilting $A$-module, $B=\End_A T$.
It is known that the correspondence 
      $$g:\dimv N\to \dimv\Hom_A(T,N)-\dimv\Ext_A^1(T,N),$$
     where $N$ is a $A$-module, induces an isomorphism $K_0(A)\to K_0(B)$ of the Grothendieck groups of $A$ and $B$. By Lemma~\ref{l:key lemma}, it is easy to get this isomorphism can be realized by $G_T^t$, \ie we have
\begin{thm}\label{t: iso induced by G_T}
Let $T$ be a tilting $A$-module, $B=\End_A T$.
Then for any $A$-module $N$, we have $$G_T^t\dimv N=\dimv\Hom_A(T,N)-\dimv\Ext_A^1(T,N).$$ 
\end{thm}
 
\begin{proof} For the sake of completeness, we provide a direct proof.

For each $1\leq j\leq n$, 
let
\begin{align}\label{gongshi: minimal projective resolution of $M$}
    0\ra P_1^j\ra P_0^j\ra T_j\ra 0
\end{align}
be the minimal projective resolution of $T_j$, where $P_1^j=\bigoplus_{i=1}^n P(i)^{b_{ij}}$, $P_0=\bigoplus_{i=1}^n P(i)^{a_{ij}}$. So  
 $g^{T_j}=[a_{1j}-b_{1j},\cdots, a_{nj}-b_{nj}]^t.$
 Then
\begin{eqnarray*}
    (g^{T_j})^t \dimv N &=&\sum_{i=1}^n(a_{ij}-b_{ij})\dim\Hom_{A}(P(i),N)\\
    &=&\sum_{i=1}^na_{ij}\dim\Hom_{A}(P(i),N)-\sum_{i=1}^nb_{ij}\dim\Hom_{A}(P(i),N)\\
    &=&\dim\Hom_{A}(P_0^j,N)-\dim \Hom_{A}(P_1^j,N).
\end{eqnarray*}
On the other hand, applying $\Hom_{A}(-,N)$ to (\ref{gongshi: minimal projective resolution of $M$}), one can get  an exact sequence
$$0\ra \Hom_{A}({T_j},N)\ra \Hom_{A}(P_0^j,N)\ra \Hom_{A}(P_1^j,N)\ra \Ext_{A}^1({T_j},N)\ra 0.$$
Then
$$\dim \Hom_{A}({T_j},N)-\dim \Ext_{A}^1({T_j},N)=\dim\Hom_{A}(P_0^j,N)-\dim \Hom_{A}(P_1^j,N).$$
So $$(g^{T_j})^t \dimv N=\dim\Hom_A({T_j},N)-\dim\Ext_{A}^1({T_j},N).$$
 Therefore
\begin{eqnarray*}
G_T^t \dimv N&=&((g^{T_j})^t\dimv N)_{1\leq j\leq n}\\
&=&(\dim\Hom_A(T_j,N)-\dim\Ext_{A}^1(T_j,N))_{1\leq j\leq n}\\
 &=&(\dim\Hom_A(T_j,N))_{1\leq j\leq n}-(\dim\Ext_{A}^1(T_j,N))_{1\leq j\leq n}\\
 &=&(\dim\Hom_B(\Hom_A(T,T_j),\Hom_A(T,N)))_{1\leq j\leq n}\\
&&-(\dim\Hom_B(\Hom_A(T,T_j),\Ext^1_A(T,N)))_{1\leq j\leq n}\\
&=&\dimv\Hom_A(T,N)-\dimv\Ext_A^1(T,N).
\end{eqnarray*} 
\end{proof}

A direct consequence of
 Theorem~\ref{t: iso induced by G_T} is
\begin{cor}\label{c: dimension vectors of $B$-modules and   $A$-modules.}
  Let  $(\mt,\mf)$ the torsion pair of $\mod A$ determined by $T$. 
  Then $$\dimv \Hom_A(T,M)=G_T^{t}\dimv M,\ \ \dimv \Ext^1_A(T,N)=-G_T^{t}\dimv N$$ for each $M\in \mt$ and $N\in\mf$.
\end{cor}





\subsection{Application to  Cartan matrices}
 \subsubsection{Congruence of Cartan matrices}
\begin{prop}\label{p: tilting module and tau-tilting module}
  Let $A$ be a finite dimensional algebra,  $T$   a tilting $A$-module, $B=\End_A T$. Then  $$\mathbf{C}_B=G_T^{t}\mathbf{C}_AG_T.$$
Moreover,  $\mathbf{det} \mathbf{C}_A=\mathbf{det} \mathbf{C}_B$. 
\end{prop}

\begin{proof}
 Let $T=\bigoplus_{i=1}^n T_i$  with each $T_i$  indecomposable. So $$\mathbf{C}_B=(\dimv\Hom_A(T,T_1),\cdots, \dimv\Hom_A(T,T_n)).$$
Because $D_T=(\dimv T_1,\cdots,\dimv T_n),$ by   Theorem~\ref{t: iso induced by G_T}, we have $$G_T^{t}D_T=(G_T^{t}\dimv T_1,\cdots, G_T^{t}\dimv T_n)=(\dimv\Hom_A(T,T_1),\cdots, \dimv\Hom_A(T,T_n))=\mathbf{C}_B.$$
 Since $\proj_A T\leq 1$, so  $D_T=\mathbf{C}_AG_{T}$. Thus, we have  $$\mathbf{C}_B=G_T^{t}D_T=G_T^{t}\mathbf{C}_AG_T.$$ Since $\mathbf{det} G_T^t\in\{1,-1\}$, we have $\mathbf{det} \mathbf{C}_A=\mathbf{det} \mathbf{C}_B$ 
 \end{proof}

 For  a finite dimensional algebra $A$ such that 
$\mathbf{det} \mathbf{C}_A\neq 0$, denote by $\Phi_A=-\mathbf{C}_A^t\mathbf{C}_A^{-1}$, $\mathbf{R}_A=(\mathbf{C}_A^{-1})^t$. 
\begin{cor}\label{c: Coxter matrxi of A and B}
     Let $T$ be a tilting $A$-module, $B=\End_A T$. Suppose that   $\mathbf{det} \mathbf{C}_A\neq 0$. Then we have
 $$\Phi_B=G_T^{t}\Phi_A(G_T^{t})^{-1} \text{ and  \   } \mathbf{R}_B=G_T^{-1}\mathbf{R}_A(G_T^{-1})^t.$$
    
\end{cor}
\begin{proof} Because $T$ is a tilting $A$-module,
by  Proposition \ref{p: tilting module and tau-tilting module}, we have $\mathbf{C}_B=G_T^{t}\mathbf{C}_AG_T$ and $\mathbf{det} \mathbf{C}_A=\mathbf{det} \mathbf{C}_B$. Because $\mathbf{det} \mathbf{C}_A\neq 0$,   $\mathbf{det} \mathbf{C}_B\neq 0$. Thus,
$$\Phi_B=-\mathbf{C}_B^t\mathbf{C}_B^{-1}=-(G_T^{t}\mathbf{C}_AG_T)^t(G_T^{t}\mathbf{C}_AG_T)^{-1}=G_T^{t}\Phi_A(G_T^{t})^{-1},$$
and   $$\mathbf{R}_B=(\mathbf{C}_B^t)^{-1}=((G_T^{t}\mathbf{C}_AG_T)^t)^{-1}=G_T^{-1}\mathbf{R}_A(G_T^{-1})^t.$$
    \end{proof}

\subsubsection{Euler form}
Assume that $A$ has finite global dimension.
Then the Euler characteristic of $A$ is the bilinear
form $\langle -,-\rangle_A$ on $K_0(A)$ defined by
$$\langle \dimv M,\dimv N\rangle_A=\sum_{s=0}^{\infty}(-1)^s\dim \Ext_A^s(M,N)$$
for any two $A$-modules $M,N$.
 It is known that $$\langle \dimv M,\dimv N\rangle_A=(\dimv M)^{t}(\mathbf{C}_A^{-1})^{t} \dimv N=(\dimv M)^{t}\mathbf{R}_A \dimv N.$$
 By Corollary~\ref{c: Coxter matrxi of A and B}, we have
\begin{cor}\label{c: roots between q_A and q_B}
    Let $A$ be a finite dimensional algebra with finite global dimension, $T$  a tilting $A$-module, $B=\End_A T$. Then for all $M,N\in\mod A$, we have
  $$\langle \dimv M,\dimv N\rangle_A=\langle G_T^{t} \dimv M, G_T^{t} \dimv N\rangle_B.$$
\end{cor}
\begin{proof}
\begin{eqnarray*}
 \langle G_T^{t} \dimv M, G_T^{t} \dimv N\rangle_B&=&(G_T^{t} \dimv M)^t\mathbf{R}_BG_T^{t} \dimv N\\
&=&(\dimv M)^{t}G_TG_T^{-1}\mathbf{R}_A(G_T^{-1})^tG_T^{t} \dimv N \\
&=&(\dimv M)^{t}\mathbf{R}_A \dimv N \\
&=&\langle \dimv M,\dimv N\rangle_A
\end{eqnarray*}
\end{proof}
\subsubsection{From $\tau$-tilting modules to tilting modules}
Let $A$ be a finite-dimensional algebra, $T$  be a $\tau$-tilting $A$-module, $B=\End_A T$. By Lemma~\ref{l:key lemma}, it is easy to get  $G_T^tD_T=\mathbf{C}_B$. Then we have

\begin{prop}\label{p: tilting module and tau-tilting module1}
    \begin{itemize}
        \item [(1)]$T$ is a tilting $A$-module if and only if  $\mathbf{C}_B=G_T^{t}\mathbf{C}_AG_T.$
        \item [(2)] Suppose that $\mathbf{det} \mathbf{C}_A\neq 0$. Then $T$ is a tilting $A$-module if and only if $\mathbf{det} \mathbf{C}_B\neq 0$ and $\mathbf{R}_B=G_T^{-1}\mathbf{R}_A(G_T^{-1})^t.$
    \end{itemize}
    
\end{prop}
\begin{proof} By Proposition~\ref{p: tilting module and tau-tilting module} and Corollary~\ref{c: Coxter matrxi of A and B},
it is enough to prove the ``if'' part.

For $(1)$, suppose $\mathbf{C}_B=G_T^{t}\mathbf{C}_AG_T$, then we have $$G_T^{t}D_T=\mathbf{C}_B=G_T^{t}\mathbf{C}_AG_T.$$ Since $\mathbf{det} G_T^t\neq 0$, then $D_T=\mathbf{C}_AG_{T}$. By Lemma~\ref{l: proj dim and g-vector},  then $\proj_A T\leq 1$. Thus $T$ is a tilting module.

For $(2)$, because $$\mathbf{C}_B=G_T^{t}\mathbf{C}_AG_T\Leftrightarrow (\mathbf{C}_B^t)^{-1}=((G_T^{t}\mathbf{C}_AG_T)^t)^{-1}\Leftrightarrow \mathbf{R}_B=G_T^{-1}\mathbf{R}_A(G_T^{-1})^t,$$
then by $(1)$, $T$ is a tilting module.
 \end{proof}

\subsection{Application to cluster-concealed algebra}
 Let \( A \) be a finite dimensional hereditary algebra, \( \mathcal{C} = D^b(\operatorname{mod} A)/\tau^{-1}[1] \) the cluster category (cf. \cite{BMRRT}). Let \( T \) be a preprojective tilting \( A \)-module, \( B = \operatorname{End}_A(T) \) a concealed algebra, and \( C = \operatorname{End}_{\mathcal{C}}(T) \) the corresponding cluster-concealed algebra. For a vector $v=(a_1,\cdots,a_n)^t$, denote by $\mathbf{abs}(v)=(|a_1|,\cdots,|a_n|)^t$.
The following is a restatement of \cite[Proposition 4 and its addendum]{R11}.
\begin{prop}
 There is a bijection $\iota: \ind(A)\to \ind(C)$ such that for each $X\in\ind(A)$, $$\dimv(\iota(X))=\mathbf{abs}(G_T^{t}\dimv X).$$
\end{prop}
\begin{proof}
By \cite[Proposition 4]{R11}, there is a bijection $\iota: \ind(A)\to \ind(C)$ such that for each $X\in\ind(A)$, $$\iota(X)|_B=\Hom_A(T,X)\oplus \Ext_A^1(T,X).$$
By \cite[Addendum to Proposition 4]{R11},
$$\dimv(\iota(X))=\mathbf{abs}(\dimv \Hom_A(T,X)-\dimv \Ext_A^1(T,X)).$$
 Then by Theorem~\ref{t: iso induced by G_T},  we have $$\dimv(\iota(X))=\mathbf{abs}(\dimv \Hom_A(T,X)-\dimv \Ext_A^1(T,X))=\mathbf{abs}(G_T^{t}\dimv X).$$
\end{proof}

\subsection{Application to mutation of certain matrices}
In this subsection, we  reformulate \cite[Corollary 7.6]{GLS06}. In particular, one can get that  the mutation on  certain kind of matrices can be realized by a $G$-matrix.
\subsubsection{Matrix mutation}Let \( B = (b_{ij}) \) be an \( l \times m \)-matrix with integer entries such that \( l \geq m \), and let \( k \in [1, m] \). Following Fomin and Zelevinsky, the {\it mutation} of \( B \) in direction \( k \) is an \( l \times m \)-matrix
$\mu_k(B) = (b'_{ij})$
defined by
\[b'_{ij} =
\begin{cases}
-b_{ij} & \text{if } i = k \text{ or } j = k, \\
b_{ij} +[-b_{ik}]_+b_{kj}+b_{ik}[b_{kj}]_+ & \text{otherwise},
\end{cases}\]
where \( 1\leq i\leq l \), \( 1\leq j\leq m \) and $[a]_+=\max\{0,a\}$. It is easy to get that $\mu_k(-B)=-\mu_k(B)$.
For an \( m \times m \)-matrix \( B \) and some \( k \in [1, m] \), we define an \( m \times m \)-matrix \( S = S(B, k) = (s_{ij}) \) by
\[s_{ij} =
\begin{cases}
-\delta_{ij} +[-b_{ij}]_+ & \text{if } i = k, \\
\delta_{ij} & \text{otherwise},
\end{cases}\]
where $\delta_{ij}$ denote the Kronecker delta.
If $B$ is a skew-symmetric \( m \times m \)-matrix,
by \cite[Lemma 7.1]{GLS06}, $S^tBS=\mu_k(B)$ for any $1\leq k\leq m$.

\subsubsection{Mutation via $G$-matrix}
 Let $\Pi$ be a classical preprojective algebra of Dynkin type, $T=T_1\oplus \cdots T_{r-n}\oplus T_{r-n+1}\oplus \cdots \oplus T_r$ be a basic maximal rigid modules with each $T_i$ indecomposable. Without loss of generality, we assume that $T_{r-n+1},\cdots,T_r$ are projective. 
 Let $E=\End_{\Pi} T$. 
 Denote by $\mathbf{C}_T$ the Cartan matrix of $E$. By \cite{GLS06}, the global dimension of $E$ is at most $3$.   Let $\mathbf{R}_T=(\mathbf{C}_{T}^{-1})^t=(r_{ij})_{r\times r}$.
By \cite[Lemma 7.3]{GLS06}, $r_{ij}=-r_{ji}$ for each $1\leq i\leq r-n$ or $1\leq j\leq r-n$. 
 Let $\Gamma_T$ be the quiver of $E$,  \( \mathbf{B}_T = (t_{ij})_{1 \leq i, j \leq r} \)  the \( r \times r \)-matrix defined by  
$$ t_{ij} = (\text{number of arrows } j \to i \text{ in } \Gamma_T) - (\text{number of arrows } i \to j \text{ in } \Gamma_T) .$$ It is clear that $\mathbf{B}_T$ is a skew-symmetric matrix. 
Let \( \mathbf{B}_T^\circ = (t_{ij}) \) and \( \mathbf{R}_T^\circ = (r_{ij}) \) be the \( r \times (r - n) \)-matrices obtained from \( \mathbf{B}_T \) and \( \mathbf{R}_T \), respectively, by deleting the last \( n \) columns. By \cite[Corollary 7.4]{GLS06}, $\mathbf{B}_T^\circ =\mathbf{R}_T^{\circ}$.

 For $1\leq k\leq r-n$, by \cite{GLS06}, there is a unique indecomposable rigid $\Pi$-module $T_k^*$ such that  $T'=T\setminus T_k\oplus T_k^*$  is also a basic maximal rigid $\Pi$-module. Let $U=\Hom_{\Pi}(T,T')$, $U$ is a tilting $E$-module and $\End_{E}U\cong \End_{\Pi}T'$ (cf. \cite{GLS06} or \cite[Theorem 5.3.2]{I07}).
Denote by  $\tilde{S}=S(-\mathbf{R}_T,k)$.  
\begin{lem}
    $G_U^t=S(-\mathbf{R}_T,k)=\tilde{S}$. 
\end{lem}
\begin{proof}
Let $\tilde{S} = (\tilde{s}_{ij})$, thus we have
\[\tilde{s}_{ij} =
\begin{cases}
-\delta_{ij} +[r_{ij}]_+ & \text{if } i = k, \\
\delta_{ij} & \text{otherwise}.
\end{cases}\]
Let $\tilde{S}^t = (\tilde{t}_{ij})$, then
\[\tilde{t}_{ij}=\tilde{s}_{ji} =
\begin{cases}
-\delta_{ji} +[r_{ji}]_+ & \text{if } j = k, \\
\delta_{ji} & \text{otherwise}.
\end{cases}\]
Because $T'=T\setminus T_k\oplus T_k^*$, there is an exact sequence $0\ra T_k\xra{g} B\xra{f} T_k^*\ra 0$, where $g$ is a minimal left $\add (T\setminus T_k)$-approximation. Then one can get a minimal projective resolution $$0\ra \Hom_{E}(T,T_k)\ra \Hom_{E}(T,B)\ra \Hom_{E}(T,T_k^*)\ra 0$$
of $\Hom_{E}(T,T_k^*)$. Because $1\leq k\leq r-n$ and $\mathbf{B}_T^\circ =\mathbf{R}_T^{\circ}$, then one can get that  $B=\bigoplus_{r_{ki}>0}T_i^{r_{ki}}$.
Thus $$g^{\Hom_{E}(T,T_k^*)}=([r_{k1}]_+,\cdots,[r_{k,k-1}]_+,-1,[r_{k,k+1}]_+,\cdots,[r_{kn}]_+)^t.$$
Because any  direct summand of $U$ except $\Hom_{E}(T,T_k^*)$ is projective, one can get that 
 $G_U=(g_{ij})$ is as following
\[g_{ij} =
\begin{cases}
-\delta_{ij} + [r_{ji}]_+ & \text{if } j = k, \\
\delta_{ij} & \text{otherwise}.
\end{cases}\]
It is clear that $G_U=\tilde{S}^t$. 
\end{proof}

By $\mathbf{B}_T^\circ =\mathbf{R}_T^{\circ}$, so for any $1\leq k\leq r-n$ and $1\leq i\leq r$, $b_{ik}=r_{ik}$. 
Because $r_{ij}=-r_{ji}$ for each $1\leq i\leq r-n$ or $1\leq j\leq r-n$ and $\mathbf{B}_T$ is skew-symmetric, 
so for $1\leq k\leq r-n$, we have $b_{kj}=-b_{jk}=-r_{jk}=r_{kj}$ for each $1\leq j\leq r$. Thus $S(-\mathbf{R}_T,k)=S(-\mathbf{B}_T,k)$. Because $\mathbf{B}_T$ is skew-symmetric, so $\tilde{S}^2=\mathbf I$, thus $G^2_U=\mathbf I$, which means $G_U^{-1}=G_U$. Now we can give another form of \cite[Proposition 7.5, Corollary 7.6]{GLS06} directly.
\begin{cor}
  $\mathbf{C}_{T'}=\tilde{S}\mathbf{C}_T\tilde{S}^t$ and $\mathbf{R}_{T'}=\tilde{S}^t\mathbf{R}_T\tilde{S}$.  
\end{cor}
\begin{proof}
By  Corollary~\ref{c: Coxter matrxi of A and B}, we have  $\mathbf{C}_{T'}=G_U^t\mathbf{C}_TG_U$ and $\mathbf{R}_{T'}=G_U^{-1}\mathbf{R}_T(G_U^{-1})^t$. Then by $G_U^t=\tilde{S}$ and $\tilde{S}^2=\mathbf I$,
 we have the results.
\end{proof}

Moreover, one can get  mutation on $\mathbf{B}_T$ can be realized by a $G$-matrix.
 \begin{thm}\label{t:mutation via g-matrix}
 Keep notation as above. Then
 $G_U\mathbf{B}_TG_U^t=\mu_k(\mathbf{B}_T)$.
 \end{thm}
\begin{proof}
By $$S(-\mathbf{B}_T,k)^t(-\mathbf{B}_T)S(-\mathbf{B}_T,k)=\mu_k(-\mathbf{B}_T),$$ 
one can get  $$S(-\mathbf{B}_T,k)^t\mathbf{B}_TS(-\mathbf{B}_T,k)=\mu_k(\mathbf{B}_T),$$
then by  $S(-\mathbf{R}_T,k)=S(-\mathbf{B}_T,k)$ and $G_U^t=S(-\mathbf{R}_T,k)$,  one can get the result.
 \end{proof}

\subsection{Application to  $g$-vectors between hereditary algebras and tilted algebras}
 Let $A$ be a finite dimensional hereditary algebra, $T$ be a tilting $A$-module, $B=\End_A T$.
In this subsection, by Theorem~\ref{t: iso induced by G_T}, we establish a relationship between $g$-vectors from $A$ and those from tilted algebra $B$.

Let $(\mt,\mf)$ the torsion pair of $\mod A$ determined by $T$.  Thus, for any $A$-module $X$, we have a unique exact sequence $0\ra tX\ra X\ra fX\ra 0$ such that $tX\in\mt$ and $fX\in \mf$. Let $(\my,\mx)$ the torsion pair of $\mod B$ determined by $T$. So $\Hom_A(T,-)$ is the equivalent functor from $\mt$ to $\my$ and  $\Ext_A^1(T,-)$ is  the equivalent functor from $\mf$ to $\mx$ (cf. \cite{ASS06}).

\begin{prop}\label{p: g-vectors between hereditary algebras and tilted algebras}
Let $X$ be an indecomposable $A$-module.
\begin{itemize} 
    \item[(1)] If $X\in\mt$, then $$g_B^{\Hom_A(T,X)}=G_T^{-1}g_A^X.$$
    \item[(2)] Let $X\in \mf$. If  $\tau_A X\in\mf$, then $$g_B^{\Ext_A^1(T,X)}=-G_T^{-1}g_A^X.$$
    \item[(3)] Let $X\in \mf$. If  $\tau_A X\notin\mf$, then $X=fP(a)$ for some indecomposable projective $A$-module $P(a)$ and  $$g_B^{\Ext_A^1(T,X)}=-G_T^{-1}g_A^{P(a)}.$$
    In particular, if $P(a)\notin \add T$, then we have $$g_B^{\Ext_A^1(T,P(a))}=-G_T^{-1}g_A^{P(a)}.$$
\end{itemize}
    
\end{prop} 
\begin{proof}
 For (1),   because  $X\in\mt$, so $\Hom_A(T,X)\in\my$, which means $\proj_B \Hom_A(T,X)\leq  1$. So $$\dimv \Hom_A(T,X)=\mathbf{C}_Bg_B^{\Hom_A(T,X)}.$$
 Note that   by Theorem~\ref{t: iso induced by G_T}, we have $\dimv \Hom_A(T,X)=G_T^{t}\dimv X$, then combined with
    $\dimv X=\mathbf{C}_Ag_A^X$,  we have 
    $$\mathbf{C}_Bg_B^{\Hom_A(T,X)}=G_T^{t}\mathbf{C}_Ag_A^X,$$
    thus, $$g_B^{\Hom_A(T,X)}=\mathbf{C}_B^{-1}G_T^{t}\mathbf{C}_Ag_A^X.$$ Associated with $\mathbf{C}_B=G_T^{t}\mathbf{C}_AG_T$, we have $$g_B^{\Hom_A(T,X)}=\mathbf{C}_B^{-1}G_T^{t}\mathbf{C}_AG_TG_T^{-1}g_A^X=G_T^{-1}g_A^X.$$ 

For (2), because $X\in \mf$ and  $\tau_A X\in\mf$,  so $\Ext_A^1(T,X)\in\mx$ and $$\tau_B\Ext_A^1(T,X)=\Ext_A^1(T,\tau_A X)\in\mx.$$ 
By Theorem~\ref{t: iso induced by G_T}, $$G_T^t\dimv\tau_AX=-\dimv\Ext_A^1(T,\tau_AX)=-\dimv \tau_B\Ext_A^1(T,X).$$
Because $\inj_A\tau_AX\leq 1$ and  $\inj_B \tau_B\Ext_A^1(T,X)\leq 1$, thus, 
by Lemma~\ref{l: proj dim and g-vector}, we have $\dimv \tau_A X=-\mathbf{C}_A^{t}g_A^X$ and  $\dimv \tau_B\Ext_A^1(T,X)=-\mathbf{C}_B^{t}g_B^{\Ext_A^1(T,X)}.$
Then we have
$$\mathbf{C}_B^{t}g_B^{\Ext_A^1(T,X)}=-G_T^{t}\mathbf{C}_A^{t}g_A^X,$$
Thus $$g_B^{\Ext_A^1(T,X)}=-(\mathbf{C}_B^{t})^{-1}G_T^{t}\mathbf{C}_A^{t}g_A^X=-(\mathbf{C}_AG_T\mathbf{C}_B^{-1})^{t}g_A^X.$$
    Associated with $\mathbf{C}_B=G_T^{t}\mathbf{C}_AG_T$, we have 
    $$g_B^{\Ext_A^1(T,X)}=-((G_T^{t})^{-1}G_T^{t}\mathbf{C}_AG_T\mathbf{C}_B^{-1})^{t}g_A^X=-G_T^{-1}g_A^X.$$

For (3), because $X\in \mf$,  $\tau_A X\notin\mf$, one can get
$\Ext_A^1(T,X)\in\mx$ but $\tau_B\Ext_A^1(T,X)\notin\mx.$ 
Then by connecting Lemma, we have $\tau_B\Ext_A^1(T,X)=\Hom_A(T,I(a))$ and  $\Ext_A^1(T,X)=\Ext_A^1(T,P(a))$ for some indecomposable projective module $P(a)$ and the corresponding indecomposable injective  module
$I(a)$.  Because $\Ext_A^1(T,fP(a))=\Ext_A^1(T,P(a))=\Ext_A^1(T,X)$,  $\Ext_A^1(T,-)$ is an equivalent functor between $\mf$ and $\mx$, so $X=fP(a)$. 
Because $T$ is a tilting $A$-module, we can take an exact sequence \begin{eqnarray}\label{g: coresolution of P(a)}
   0\ra P(a)\ra T_0\ra T_1\ra 0. 
\end{eqnarray}
 Moreover, since $X=fP(a)\neq 0$, which means $P(a)\notin \add T$, so $P(a)\notin \add T_0$. Suppose $tP(a)\neq 0$, because $A$ is a  hereditary algebra, so $tP(a)$ is a projective $A$-module. Thus, $tP(a)\in \add T$.
Applying $\Hom_A(T,-)$ to (\ref{g: coresolution of P(a)}), one can get an exact sequence 
\begin{eqnarray}\label{gongshi: projective resolution of ext p}
    0\ra \Hom_A(T,P(a))\ra \Hom_A(T,T_0)\ra \Hom_A(T,T_1)\ra \Ext_A^1(T,P(a))\ra 0,
\end{eqnarray}
which is a projective resolution of $\Ext_A^1(T,P(a))$ since $\Hom_A(T,P(a))=\Hom_A(T,tP(a))$.
Moreover, because $A$ is a finite dimensional  hereditary algebra, so no direct summand of $tP(a)$ lies in $\add T_0$, thus no direct summand of 
$\Hom_A(T,P(a))$—or equivalently, no direct summand of $\Hom_A(T,tP(a))$— lies in $\add \Hom_A(T,T_0)$. Therefore, one can get   $$\mathbf{C}_Bg_B^{\Ext_A^1(T,P(a))}=\dimv \Hom_A(T,T_1)- \dimv \Hom_A(T,T_0).$$
Then by Theorem~\ref{t: iso induced by G_T} and (\ref{gongshi: projective resolution of ext p}), we have
\begin{eqnarray*}
   G_T^{t}\dimv P(a)&=&\dimv \Hom_A(T,P(a))-\dimv \Ext_A^1(T,P(a)) \\
   &=&\dimv \Hom_A(T,T_0)- \dimv \Hom_A(T,T_1)\\
   &=&-\mathbf{C}_Bg_B^{\Ext_A^1(T,P(a))}.
\end{eqnarray*}
Hence, $$G_T^{t}\mathbf{C}_Ag_A^{P(a)}=G_T^{t}\dimv P(a)=-\mathbf{C}_Bg_B^{\Ext_A^1(T,P(a))},$$
which means $$g_B^{\Ext_A^1(T,P(a))}=-\mathbf{C}_B^{-1}G_T^{t}\mathbf{C}_Ag_A^{P(a)}.$$
By $\mathbf{C}_B=G_T^{t}\mathbf{C}_AG_T$,  one can get
$$g_B^{\Ext_A^1(T,X)}=g_B^{\Ext_A^1(T,fP(a))}=g_B^{\Ext_A^1(T,P(a))}=-G_T^{-1}g_A^{P(a)}.$$


\end{proof}

By Proposition~\ref{p: g-vectors between hereditary algebras and tilted algebras}, we have 
\begin{cor}
For $1\leq i\leq n$, denote 
$g_i=
\begin{cases}
g_B^{\Hom_A(T,P(i))} & \text{if } P(i)\in\add T, \\
-g_B^{\Ext^1_A(T,P(i))} & \text{otherwise}.
\end{cases}$
Then $$G_T^{-1}=(g_1,\cdots,g_n).$$
\end{cor}
\begin{proof}
Without loss of generality, we can suppose that there is an integer $j\in [1,n]$ such that $P(i)\in\add T$    for $1\leq i\leq j$ and $P(i)\notin \add T$ for $j<i\leq n$.
Then we have 
\begin{eqnarray*}
(g_1,\cdots,g_n)&=&(g_B^{\Hom_A(T,P(1))},\cdots,g_B^{\Hom_A(T,P(j))},- g_B^{\Ext^1_A(T,P(j+1))},\cdots, -g_B^{\Ext^1_A(T,P(n))})\\
&=&(G_T^{-1}g_A^{P(1)},\cdots,G_T^{-1}g_A^{P(j)}, G_T^{-1}g_A^{P(j+1)},\cdots,G_T^{-1}g_A^{P(n)})\\
&=&G_T^{-1}(g_A^{P(1)},\cdots,g_A^{P(n)})\\
&=&G_T^{-1}.
\end{eqnarray*}

\end{proof}

\section{G-Matrix views on  Coxeter, Nakayama and AR-Translation}\label{s:other functor-induced iso and G-matrix}
 In this section, for a finite dimensional  algebra $A$, we fix $G_{DA}=(g^{I(1)},\cdots,g^{I(n)})$ and we will show that for some special algebras,  the Coxeter matrix, the Nakayama functor and the Auslander-Reiten functor can also be realized by $G$-matrices.
\subsection{Coxeter transformations via G-matrices for 1-Gorenstein algebras}
Let $A$ be a finite dimensional  algebra.  $A$ is called {\it $1$-Gorenstein} if $\proj DA\leq 1$.
\begin{prop}\label{p:judge 1-gorenstien alg}
The following are equivalent
 \begin{itemize}
     \item[(1)]  $A$ is $1$-Gorenstein algebra;
     \item[(2)] $\mathbf{C}_AG_{DA}=\mathbf{C}_A^{t}$;
      \item[(3)] $DA$ is a tilting $A$-module;
     \item[(4)] $DA$ is a $\tau$-tilting $A$-module.
    
\end{itemize}
 Moreover, if $\mathbf{det} \mathbf{C}_A\neq 0$, then $(1),(2)$, $(3)$, $(4)$ are also equivalent to
 \begin{itemize}
     \item [(5)] $\Phi_A=-(G_{DA}^{-1})^{t}$.
 \end{itemize}
In particular,  in this case,   $\Phi_A$ has column  sign-coherence property.
\end{prop}
\begin{proof}
By  Lemma~\ref{l: proj dim and g-vector}, $\proj_A DA\leq 1$ if and only if $\mathbf{C}_AG_{DA}=\mathbf{C}_A^{t}.$
So $(1),(2)$ and $(3)$ are equivalent. It is clear that if $DA$ is a tilting $A$-module, then $DA$ is a $\tau$-tilting $A$-module.
 On the other hand, if $DA$ is $\tau$-tilting, because $DA$ is faithful, so $DA$ is tilting. Thus  $(1),(2)$, $(3)$, $(4)$ are equivalent.

Now assume that $\mathbf{det} \mathbf{C}_A\neq 0$. 

\noindent $(2)\Rightarrow(5)$, by $\mathbf{C}_AG_{DA}=\mathbf{C}_A^{t}$, we have $G_{DA}=\mathbf{C}_A^{-1}\mathbf{C}_A^{t}$, then $$(G_{DA}^{-1})^{t}=((\mathbf{C}_A^{-1}\mathbf{C}_A^{t})^{-1})^t=\mathbf{C}_A^{t}\mathbf{C}_A^{-1}=-\Phi_A.$$
\noindent $(5)\Rightarrow(2)$, by $(G_{DA}^{-1})^{t}=-\Phi_A=\mathbf{C}_A^{t}\mathbf{C}_A^{-1},$ one can get that $G_{DA}=\mathbf{C}_A^{-1}\mathbf{C}_A^{t}$.
Therefore, $(1),(2)$, $(3)$, $(4)$ and $(5)$ are equivalent in this case.
Because $DA$ is a tilting module, $G_{DA}$ is a $G$-matrix, so  $\Phi_A=-((G_{DA})^{-1})^{t}$ is a 
$C$-matrix, thus $\Phi_A$ has column sign-coherence property.
\end{proof}

\subsection{Nakayama functor via $G$-matrix for self-injective algebra}
Let $A$ be   a self-injective algebra, $\nu(-)=D\Hom_A(-,A)$,  the Nakayama functor, $\nu^{-1}=\Hom_A(DA,-)$ the quasi-inverse of $\nu$.
Note that there is a Nakayama permutation $\sigma$: $Q_0\ra Q_0$ of $A$, \ie $e_iA=D(A{e_{\sigma(i)}})$. 
So $g^{I(i)}=\mathbf e_{\sigma^{-1}(i)}$  for each $i$.  Then we have
\begin{lem}
Let $A=kQ/I $ be a finite dimensional algebra. Then $A$ is a self-injective algebra if and only if there is a permutation $\sigma$: $Q_0\ra Q_0$ such that 
$$G_{DA}=(\mathbf e_{\sigma^{-1}(1)},\cdots,\mathbf e_{\sigma^{-1}(n)}).$$ 
\end{lem}
\begin{proof}
Suppose that $A$ is a self-injective algebra,  so  there is a Nakayama permutation $\sigma: Q_0\ra Q_0$ of $A$, such that  $e_iA=D(A{e_{\sigma(i)}})$. 
Thus $g^{I(i)}=g^{D(A{e_{i}})}=g^{e_{\sigma^{-1}(i)}A}=\mathbf e_{\sigma^{-1}(i)}$. Hence, $$G_{DA}=(\mathbf e_{\sigma^{-1}(1)},\cdots,\mathbf e_{\sigma^{-1}(n)}).$$
On the other hand, if there  is a permutation $\sigma$: $Q_0\ra Q_0$ such that 
$$G_{DA}=(\mathbf e_{\sigma^{-1}(1)},\cdots,\mathbf e_{\sigma^{-1}(n)}), $$
So each column of $G_{DA}$ is non-negative, which means that each injective $A$-module is a projective $A$-module, so $A$ is a self-injective algebra. 
\end{proof}

\begin{prop}\label{p: nakayama for self-injective algebra}
Let $A$ be a self-injective algebra. For any   $X\in \mod A$, we have
\begin{itemize}
    \item[(1)] $\dimv \nu X=G_{DA}\dimv X=(G_{DA}^{-1})^{t}\dimv X$;
    \item[(2)] $ g^{\nu X}=G_{DA}g^X=(G_{DA}^{-1})^{t}g^X.$
\end{itemize}
\end{prop}
\begin{proof}
For $(1)$, suppose that $\dimv X=(c_1,\cdots,c_n)^t$, $\dimv \nu X=(d_1,\cdots,d_n)^t$. Because $\nu$ is functorially isomorphic to $-\otimes_ADA$, so $$d_i=\dim\Hom_A(\nu X,I(i))=\dim\Hom_A(X\otimes_ADA,I(i))=\dim\Hom_A(X,\Hom_A(DA,I(i))).$$
Because $\Hom_A(DA,I(i))=\nu^{-1}I(i)=P(i)$, so
$$d_i=\dim\Hom_A(X,\Hom_A(DA,I(i)))=\dim\Hom_A(X,P(i))=\dim\Hom_A(X,I(\sigma(i)))=c_{\sigma(i)}.$$
Note that $$(\mathbf e_{\sigma^{-1}(1)},\cdots,\mathbf e_{\sigma^{-1}(n)})(c_1,\cdots,c_n)^t=c_1\mathbf e_{\sigma^{-1}(1)}+\cdots+c_n\mathbf e_{\sigma^{-1}(n)}=(c_{\sigma(1)},\cdots,c_{\sigma(n)})^t,$$
$G_{DA}=(\mathbf e_{\sigma^{-1}(1)},\cdots,\mathbf e_{\sigma^{-1}(n)}),$
so we have $\dimv \nu X=G_{DA}\dimv X$. Since $G_{DA}$ is a permutation matrix, so $G_{DA}$ is an orthogonal matrix, then we have $$\dimv \nu X=G_{DA}\dimv X=(G_{DA}^t)^{-1}\dimv X.$$

For $(2)$, let $P_1\xra{f} P_0\ra X\ra 0$ be a minimal projective resolution of $X$, where $P_0=\bigoplus_{i=1}^nP(i)^{a_i}$ and  $P_1=\bigoplus_{i=1}^nP(i)^{b_i}$. So $g^X=(a_1-b_1,\cdots,a_n-b_n)^t$ and 
\begin{eqnarray*}
    G_{DA}g^X&=&(\mathbf e_{\sigma^{-1}(1)},\cdots,\mathbf e_{\sigma^{-1}(n)})(a_1-b_1,\cdots,a_n-b_n)^t\\
    &=&(a_1-b_1)\mathbf e_{\sigma^{-1}(1)}+\cdots+(a_n-b_n)\mathbf e_{\sigma^{-1}(n)}\\
    &=&(a_{\sigma(1)}-b_{\sigma(1)},\cdots,a_{\sigma(n)}-b_{\sigma(n)})^t.
\end{eqnarray*}
Note that $$\nu P_1\xra{\nu f} \nu P_0\ra \nu X\ra 0$$ is an minimal projective resolution  of $\nu X$. Because $\nu P(i)=I(i)=D(A{e_{i}})=e_{\sigma^{-1}(i)}A=P(\sigma^{-1}(i))$, So $$\nu P_0=\bigoplus_{i=1}^nP(\sigma^{-1}(i))^{a_i}=\bigoplus_{\sigma^{-1}(i)=1}^nP(i)^{a_{\sigma(i)}},\ \ \nu P_1=\bigoplus_{i=1}^nP(\sigma^{-1}(i))^{b_i}=\bigoplus_{\sigma^{-1}(i)=1}^nP(i)^{b_{\sigma(i)}},$$
thus we have $g^{\nu X}=(a_{\sigma(1)}-b_{\sigma(1)},\cdots,a_{\sigma(n)}-b_{\sigma(n)})^t=G_{DA}g^X$.  Then using $G_{DA}$ is an orthogonal matrix,   $(2)$ is true.
\end{proof}

\begin{rk} Let $A$ be a finite-dimensional self-injective algebra, 
if  $\mathbf{det}\mathbf{C}_A\neq 0$, by Proposition~\ref{p:judge 1-gorenstien alg}, for any $X\in\mod A$, one can get that 
$$\Phi_A(\dimv X)=-(G_{DA}^t)^{-1}\dimv X=-\dimv\nu X.$$ 
\end{rk}

\subsection{AR-translation via $G$-matrix for  hereditary algebras}
 \begin{prop}\label{p: tau induced by g-vector for hereditary alg}
    Let $A$ be a finite dimensional hereditary algebra. 
    \begin{itemize}
    \item [(1)] For any $A$-module $M$, $\dimv M=\mathbf{C}_Ag^M$.
        \item [(2)] For any indecomposable non-projective $A$-module $M$,\begin{itemize}
            \item [(2.1)]$\dimv \tau M=-\mathbf{C}_A^{t}g^M=-((G_{DA})^{-1})^{t}(\dimv M);$
            \item [(2.2)]$g^{\tau M}=-G_{DA}(g^M).$
        \end{itemize}
         \item [(3)] For any indecomposable non-injective $A$-module $M$,\begin{itemize}
            \item [(3.1)]$\dimv \tau^{-1} M=-(\mathbf{C}_A^{t})^{-1}g^M=-G_{DA}^{t}(\dimv M);$
            \item [(3.2)]$g^{\tauni M}=-G_{DA}^{-1}(g^M).$
        \end{itemize}
    \end{itemize}
    \end{prop}
  \begin{proof}
   Because $A$ is hereditary, so for any $A$-module $M$, $\proj_A M\leq 1$, thus $\dimv M=\mathbf{C}_Ag^M$. In particular, $\mathbf{C}_AG_{DA}=\mathbf{C}_A^t$, thus one can get  $G_{DA}=\mathbf{C}_A^{-1}\mathbf{C}_A^t$.
   
Because $(1)$ is clear, $(3)$ can be get from $(2)$ directly, so we only need to prove $(2)$.
   
  For (2),  because  $A$ is hereditary, so  $\inj_A\tau M\leq 1$. Then by Lemma~\ref{l: proj dim and g-vector},
$\dimv \tau M=-\mathbf{C}_A^{t}g^M.$
By $G_{DA}=\mathbf{C}_A^{-1}\mathbf{C}_A^t$,  one can get $$((G_{DA})^{-1})^{t}(\dimv M)=((\mathbf{C}_A^{-1}\mathbf{C}_A^t)^{-1})^{t}(\dimv M)=\mathbf{C}_A^{t}\mathbf{C}_A^{-1}(\dimv M)=\mathbf{C}_A^{t}g^M.$$ So $(2.1)$ is true.
By $\dimv \tau M=\mathbf{C}_Ag^{\tau M}$, we have
$$g^{\tau M}=\mathbf{C}_A^{-1}\dimv \tau M=-\mathbf{C}_A^{-1}\mathbf{C}_A^{t}g^M=-G_{DA}(g^M).$$
Hence, $(2.2)$ is true.
\end{proof}

\section{Geometric representations of symmetric groups via G-matrices}\label{s:Symmetric groups via G-matrices}

\subsection{Symmetric groups}\label{s:def of symmetric group}
Now we recall some well-known properties of the symmetric groups. We consider the action of \(\mathfrak{S}_{n+1}\) on \(\mathbb{R}^{n+1}\) given by permuting the standard basis \(\mathbf e_1, \ldots, \mathbf e_{n+1}\). Then \(\mathfrak{S}_{n+1}\) acts on the subspace
\[V := \left\{ x_1\mathbf e_1 + \cdots + x_{n+1}\mathbf e_{n+1} \in \mathbb{R}^{n+1} \ \middle| \ \sum_{i=1}^n x_i = 0 \right\},\]
which has a basis \(\alpha_i :=\mathbf e_i -\mathbf e_{i+1}\) with \(1 \leq i \leq n\). Clearly the action of \(\mathfrak{S}_{n+1}\) on \(V\) is faithful, and we have an injective homomorphism \(\mathfrak{S}_{n+1} \to \mathrm{GL}(V)\) called the {\it geometric representation}.

Let \(s_i\) be the transposition \((i, i+1) \in \mathfrak{S}_{n+1}\). Then \[s_{i}(\alpha_j) =
\begin{cases}
-\alpha_i & \text{if } j=i, \\
\alpha_{i}+\alpha_j & \text{else if } |j-i|=1,\\
\alpha_j & \text{otherwise.} 
\end{cases}\]
Hence, we have 
 \begin{eqnarray*}
  s_i(\alpha_1,\cdots,\alpha_n)&= &(\alpha_1,\cdots,\alpha_n)\begin{bmatrix}
1 &  &  &  & &&\\
 & \ddots &   & & &\\
  &  &  1& & &  & \\
 && 1& -1 &1& & \\
 &  &  &  &1&& \\
  &  &  &  &&\ddots& \\
   &  &  &  &&& 1
\end{bmatrix}\\
&=&(\alpha_1,\cdots,\alpha_n)R_i,  
\end{eqnarray*}
where $-1$ is in the $i$-th column of  $R_i$. 
For a  word $w=s_{i_1}\cdots s_{i_l}=s_{j_1}\cdots s_{j_m}\in\mathfrak{S}_{n+1}$, it is easy to get that $$R_{i_1}\cdots R_{i_l}=R_{j_1}\cdots R_{j_m}.$$
Then we can  set $$R_w=R_{i_1}\cdots R_{i_l},$$ which is  called the {\it geometric representation matrix} of  $w$. Then one can get $$w(\alpha_1,\cdots,\alpha_n)=s_{i_1}\cdots s_{i_l}(\alpha_1,\cdots,\alpha_n)=(\alpha_1,\cdots,\alpha_n)R_{i_1}\cdots R_{i_l}=(\alpha_1,\cdots,\alpha_n)R_w.$$
 It is clear that \(\mathfrak{S}_{n}\) is isomorphic to the subgroup of \(\mathfrak{S}_{n+1}\) generated by \(s_1, s_2, \ldots, s_{n-1}\). In the following, we identify \(\mathfrak{S}_{n}\) with this subgroup of \(\mathfrak{S}_{n+1}\). That is, we regard \(\mathfrak{S}_{n}\) as the subgroup of \(\mathfrak{S}_{n+1}\) consisting of all permutations that fix the element \(n+1\).

\subsection{Tilting modules on  Auslander algebras of $k[x]/(x^n)$}
Let $\Lambda$ be the Auslander algebras of $k[x]/(x^n)$.
\ie $\Lambda=kQ/I$, where $Q$ is as follows,
\[\xymatrix{
1\ar@<.5ex>[r]^{a_1}&2\ar@<.5ex>[r]^{a_2}\ar@<.5ex>[l]^{b_2}&3\ar@<.5ex>[l]^{b_3}\ar@<.5ex>[r]&\cdots\ar@<.5ex>[l]\ar@<.5ex>[r]&n-1\ar@<.5ex>[r]^{a_{n-1}}\ar@<.5ex>[l]& n\ar@<.5ex>[l]^{b_n}
}
\]
and $I$ is generated by $a_1b_2=0$ and $a_ib_{i+1}=b_{i}a_{i-1}$ for each $2\leq i\leq n-1$.
For each $1\leq i\leq n$, denoted by $I_i=\Lambda (1-e_i)\Lambda$ the two-sided ideal. Denoted by $\langle I_1,\cdots,I_{n-1}\rangle$ the ideal semigroup generated by $I_1,\cdots,I_{n-1}$.
\begin{lem}\cite{IZ20}
    Let $\Lambda$ be the Auslander algebra of $k[x]/(x^n)$.
\begin{enumerate}
    \item The set of tilting $\Lambda$-modules is given by $\langle I_1, \dots, I_{n-1} \rangle$.
    \item There exists a well-defined bijection $I: \mathfrak{S}_n \to \langle I_1, \dots, I_{n-1} \rangle$, which maps $w$ to
    \[
    I_w = I_{i_1} \cdots I_{i_t},
    \]
    where $w = s_{i_1} \cdots s_{i_t}$ is an arbitrary reduced expression.
    Moreover, $\End_{\Lambda} I_w=\Lambda$.
    \end{enumerate}
    \end{lem}

\begin{lem}\label{l: resolusion of e_iI_j}\cite{IZ20}
     Let \(\Lambda\) be the Auslander algebra of \( K[x]/(x^n) \). 

\begin{enumerate}
  
\item[(1)] For \( 1 \leq i \leq n - 1 \), there exist minimal projective resolution of $\mathrm{rad} P(i)$
    \[
    0 \to P(i) \to P(i-1) \oplus P(i+1) \to \operatorname{rad} P(i) \to 0.
    \]

    \item[(2)] There exist minimal projective resolution of $\mathrm{rad} P(n)$
    \[
    0 \to P({n-1}) \to \operatorname{rad} P(n) \to 0,
    \]
\end{enumerate}
where $\operatorname{rad} P(i)$ means the radical of $P(i)$.
\end{lem}

\subsection{Geometric representations via G-matrices}
The Grothendieck group \(K_0(\Lambda)\) is a free abelian group with basis \([S_1], \ldots, [S_{n}]\). 
Since \( I_w \) is a tilting \(\Lambda\)-module with \(\operatorname{End}_\Lambda(I_w) \cong \Lambda\) for any \( w \in \mathfrak{S}_n \), we have an autoequivalence
\[
\mathbf{R}\Hom_{\Lambda}(I_w,-) : D^b(\operatorname{mod} \Lambda) \to D^b(\operatorname{mod} \Lambda)
\]
which induces an automorphism $[\mathbf{R}\Hom_{\Lambda}(I_w,-)]$ on \(K_0(\Lambda)\).
For any $w\in \mathfrak{S}_n$, we fix $$G_{I_w}=(g^{e_1I_w},\cdots,g^{e_nI_w}).$$
Because  \[
I_i=e_1I_i\oplus\cdots \oplus e_nI_i= P(1) \oplus \cdots \oplus P(i-1)\oplus \operatorname{rad} P(i)\oplus P(i+1) \oplus \cdots \oplus P(n),
\]
By Lemma~\ref{l: resolusion of e_iI_j}, it is easy to get that 
\begin{lem} \label{l: G and R for i}
For $1\leq i\leq n-1$, 
    $G_{I_i}^t=R_i$.
\end{lem}
Moreover, we have
\begin{thm}\label{t: symmetric group via g-matrices}
    For any word  $w\in \mathfrak{S}_n$, we have $G_{I_w}^t=R_w$.
\end{thm}
\begin{proof}
Note that for any  $M\in\mod \Lambda$, by $\proj_{\Lambda}I_w\leq 1$ and Theorem~\ref{t: iso induced by G_T}, we have
$$[\mathbf{R}\Hom_{\Lambda}(I_w,M)]=\dimv\Hom_{\Lambda}(I_w,M)-\dimv\Ext_{\Lambda}^1(I_w,M)=G_{I_w}^t\dimv M.$$
So $$[\mathbf{R}\Hom_{\Lambda}(I_w,-)]([M])=[\mathbf{R}\Hom_{\Lambda}(I_w,M)]=G_{I_w}^t\dimv M.$$
Thus \begin{eqnarray*}
    [\mathbf{R}\Hom_{\Lambda}(I_w,-)]([S_1],\cdots,[S_n])&=&(G_{I_w}^t\dimv S_1, \cdots, G_{I_w}^t\dimv S_n)\\
    &=&G_{I_w}^t(\dimv S_1, \cdots, \dimv S_n)\\
    &=&G_{I_w}^t.
\end{eqnarray*}
For $1\leq i\leq n-1$, by Lemma~\ref{l: G and R for i}, we have
 $$[\mathbf{R}\Hom_{\Lambda}(I_i,-)][S_j]=[\mathbf{R}\Hom_{\Lambda}(I_i,S_j)]=G_{I_i}^t\dimv S_j=R_i\dimv S_j.$$
Hence, $$[\mathbf{R}\Hom_{\Lambda}(I_i,-)]([S_1],\cdots,[S_n])=R_i(\dimv S_1,\cdots,\dimv S_n)=R_i.$$
Note that for any reduced expression  $w = s_{i_1}s_{i_2}\cdots s_{i_k}$, $I_w = I_{i_1}I_{i_2}\cdots I_{i_k}=I_{i_1}\lten_{\Lambda}I_{i_2}\lten_{\Lambda}\cdots \lten_{\Lambda} I_{i_k}$ (cf. \cite{IZ20}).
So \begin{eqnarray*}
    \mathbf{R}\Hom_{\Lambda}(I_w,-)&=&\mathbf{R}\Hom_{\Lambda}(I_{i_1}\lten_{\Lambda}I_{i_2}\lten_{\Lambda}\cdots \lten_{\Lambda} I_{i_k},-)\\
    &=&\mathbf{R}\Hom_{\Lambda}(I_{i_1},\cdots, \mathbf{R}\Hom_{\Lambda}(I_{i_{k-1}}, \mathbf{R}\Hom_{\Lambda}(I_{i_k},-))\cdots).
\end{eqnarray*}
Then 
\begin{eqnarray*}
    [\mathbf{R}\Hom_{\Lambda}(I_w,-)]&=&[\mathbf{R}\Hom_{\Lambda}(I_{i_1},-)]\cdots[\mathbf{R}\Hom_{\Lambda}(I_{i_k},-)].
\end{eqnarray*}
So \begin{eqnarray*}
    [\mathbf{R}\Hom_{\Lambda}(I_w,-)]([S_1],\cdots,[S_n])&=&[\mathbf{R}\Hom_{\Lambda}(I_{i_1},-)]\cdots[\mathbf{R}\Hom_{\Lambda}(I_{i_k},-)]([S_1],\cdots,[S_n])\\
    &=&R_{i_1}\cdots R_{i_k}\\
    &=&R_{w}.    
\end{eqnarray*}
Therefore, $G_{I_w}^t=R_w$. 
    \end{proof}
\begin{cor}
   For any word $w\in\mathfrak{S}_n$,  $R_{w}$ has column sign-coherence.
\end{cor}   

For any  two words  $w, w'\in \mathfrak{S}_n$, $I_{ww'}$ and  $I_wI_{w'}$ are not equal in general, but $R_{{ww'}}=R_wR_{w'}$ always holds, then one can get
\begin{cor}\label{c: multiplycation of G-matrix}
For any  two words  $w, w'\in \mathfrak{S}_n$, then 
  $G_{I_{ww'}}=G_{I_{w'}}G_{I_w}$.
\end{cor}
\begin{proof}
  By Theorem~\ref{t: symmetric group via g-matrices},  $G_{I_{ww'}}=R_{{ww'}}^t=(R_wR_{w'})^t=(R_{w'})^tR_w^t=G_{I_{w'}}G_{I_w}$.
\end{proof}



\section{Geometric representations of Weyl group via G-matrices}\label{s:Weyl groups via G-matrices}
In this section, for a Kac-Moody Lie algebra,  we show that each reflection can be realized by a $G$-matrix from the generalized preprojective algebra.

\subsection{Weyl group and reflection}\label{ss:Weyl group}
Let $C = (c_{ij}) \in M_n(\mathbb{Z})$ be a {\it symmetrizable (generalized)
Cartan matrix} with a symmetrizer  $D = \operatorname{diag}(c_1,\cdots , c_n)$ with $c_i \geq 1$ for all $i$.
Let $\Phi=\Phi(C)$ be the root system of $C$. Let $\{\alpha_1,\cdots,\alpha_n\}\subset \Phi$ be a set of simple roots and $L$ the root lattice.  We let $V=L\otimes_{\mathbb{Z}}\mathbb{R}$ and denote by $V^*$ the dual of $V$ with the basis  $\{\alpha_1^*,\cdots,\alpha_n^*\}$.
We denote the natural pairing by  $(v^*,v)$  for $v\in V$ and $v^*\in V^*$. Then for any $1\leq i,j \leq n$, we have $(\alpha_j^*,\alpha_i)=\delta_{ij}$.

For each $1\leq i\leq n$,
define a reflection $s_i: V\ra V$ by $$s_i(\alpha_j):=\alpha_j-c_{ij}\alpha_i$$ for any $1\leq j\leq n$.
Then one can get that \begin{eqnarray*}
  s_i(\alpha_1,\cdots,\alpha_n)&= &(\alpha_1,\cdots,\alpha_n)\begin{bmatrix}
1 &  &  &  & &&\\
 & \ddots &   & & &\\
  &  &  1& & &  & \\
-c_{i1} &\cdots & -c_{i,i-1}& -1 &-c_{i,i+1}&\cdots & -c_{in} \\
 &  &  &  &1&& \\
  &  &  &  &&\ddots& \\
   &  &  &  &&& 1
\end{bmatrix}\\
&=&(\alpha_1,\cdots,\alpha_n)R_i  
\end{eqnarray*}
The Weyl group is defined as a subgroup
 $W=W(C)=\langle s_1,\cdots,s_n\rangle$ of $\mathrm{GL}(V)$.
 
 For each word $w\in W(C)$, define $\rho(w)\in \mathrm{GL}{(V^*)}$ such that $$(\rho(w)(f))(v)=(f,w^{-1}v), \text{\ for } f\in V^*, v\in V,$$
 which means $(\rho(w)(f),v)=(f,w^{-1}v)$ for any $w\in W(C), f\in V^*, v\in V$.
Let $w_1,w_2\in W(C)$, then for any $f\in V^*, v\in V$, we have $$(\rho(w_1w_2)(f))(v)=(f,w_2^{-1}w_1^{-1}v)$$ and 
$$(\rho(w_1)\rho(w_2)(f))(v)=(\rho(w_1)(\rho(w_2)(f)))(v)=(\rho(w_2)(f),w_1^{-1}v)=(f,w_2^{-1}w_1^{-1}v).$$
Then one can get that $\rho(w_1w_2)=\rho(w_1)\rho(w_2)$ for any $w_1,w_2\in W(C)$.
Hence, for any word $w=s_{i_1}\cdots s_{i_k}$, we have $$\rho(w)=\rho(s_{i_1}\cdots s_{i_k})=\rho(s_{i_1})\cdots \rho(s_{i_k}).$$
Let $p=\sum\limits_{j=1}^np_j\alpha_j^*\in V^*$. Because for each $1\leq k\leq n$, $$\rho(s_i)(p)(\alpha_k)=(p,s_i(\alpha_k))=(p,\alpha_k-c_{ik}\alpha_i)=(\sum\limits_{j=1}^np_j\alpha_j^*, \alpha_k-c_{ik}\alpha_i)=p_k-p_ic_{ik}.$$
Then we have $\rho(s_i)(p)=\sum\limits_{k=1}^n(p_k-p_ic_{ik})\alpha_k^*=p-p_i(\sum\limits_{k=1}^nc_{ik}\alpha_k^*)$ for  $p=\sum\limits_{j=1}^np_j\alpha_j^*$.
In particular, we have
\[\rho(s_i)(\alpha_l^*)=\begin{cases}\alpha_i^*-\sum\limits_{j=1}^n c_{ij}\alpha_j^* & if\ l=i;  \\ \alpha_l^* & if\ l\neq i,\end{cases}\]
which means
\begin{eqnarray*}
      \rho(s_i)(\alpha_1^*,\cdots,\alpha_n^*)&=&(\alpha_1^*,\cdots,\alpha_n^*)\begin{bmatrix}
1 &  &  &-c_{i1}  & &&\\
 & \ddots & &\vdots  & & &\\
  &  &  1& -c_{i,i-1}& &  & \\
 & & & -1 && &  \\
 &  &  &  -c_{i,i+1} &1&& \\
  &  &  & \vdots&&\ddots& \\
   &  &  & -c_{in} &&& 1
\end{bmatrix}\\
&=&(\alpha_1^*,\cdots,\alpha_n^*)\Sigma_i  
\end{eqnarray*}
It is clear that $\Sigma_i=R_i^{t}.$

For a word \(w=s_{i_{1}}\cdots s_{i_{k}}=s_{j_1}\cdots s_{j_m}\in W(C)\), 
it is easy to get that $$R_{i_1}\cdots R_{i_k}=R_{j_1}\cdots R_{j_m}.$$
Then we can
set  
\begin{itemize}
    \item \(R_{w}=R_{i_{1}}\cdots R_{i_{k}}\), called the {\it geometric representation matrix} of \(w\);  
    \item \(\varSigma_{w}=\varSigma_{i_{1}}\cdots \varSigma_{i_{k}}\), called the {\it  contragredient matrix} of \(w\).
\end{itemize}
 Then one can get $$w(\alpha_1,\cdots,\alpha_n)=s_{i_1}\cdots s_{i_k}(\alpha_1,\cdots,\alpha_n)=(\alpha_1,\cdots,\alpha_n)R_{i_1}\cdots R_{i_k}=(\alpha_1,\cdots,\alpha_n)R_w$$ and
$$ \rho(w)(\alpha_1^*,\cdots,\alpha_n^*)=\rho(s_{i_1})\cdots\rho(s_{i_k})(\alpha_1^*,\cdots,\alpha_n^*)=(\alpha_1^*,\cdots,\alpha_n^*)\Sigma_{i_1}\cdots \Sigma_{i_k}=(\alpha_1^*,\cdots,\alpha_n^*)\Sigma_w.$$ 
Moreover, we have $\Sigma_w=\Sigma_{i_1}\cdots \Sigma_{i_k}=R^t_{i_1}\cdots R^t_{i_k}=(R_{i_k}\cdots R_{i_1})^t=R_{w^{-1}}^t$.

\subsection{Generalized preprojective algebra}
Let $C = (c_{ij}) \in M_n(\mathbb{Z})$ be a {\it symmetrizable (generalized)
Cartan matrix} with a symmetrizer  $D = \operatorname{diag}(c_1,\cdots , c_n)$ with $c_i \geq 1$ for all $i$.
Denote
$g_{ij} := | \operatorname{gcd}(c_{ij} , c_{ji})|, f_{ij} := |c_{ij} |/g_{ij} .$

An {\it orientation} of $C$ is a subset of
$\{1, 2,\cdots , n\} \times \{1, 2, \cdots , n\} $ such that the followings hold:

$(A1)$~ $\{(i, j), (j, i)\} \cap \Omega
 \neq \emptyset $ if and only if $c_{ij} < 0$;

$(A2)$~ For each sequence $((i_1, i_2), (i_2, i_3), \cdots , (i_t, i_{t+1}))$ with $t \geq 1$ and $(i_s, i_{s+1}) \in \Omega$
 for
all $1 \leq s \leq t$, we have $i_1 \neq i_{t+1}$.

Given an orientation $\Omega$
 of $C$, let $Q := Q(C,\Omega
) := (Q_0,Q_1, s, t)$ be the quiver with the set of
vertices $Q_0 := \{1, \cdots , n\}$ and with the set of arrows
$Q_1 := \{a^{(g)}_{ij} : j \to i | (i, j) \in \Omega, 1 \leq g \leq g_{ij}\} \cup\{\varepsilon_i : i \to i | 1 \leq i \leq n\}.$
Thus we have $s(a^{(g)}_{ij} ) = j$ and $t(a^{(g)}_{ij} ) = i$ and $s(\varepsilon_i) = t(\varepsilon_i) = i$, where $s(a)$ and $t(a)$
denote the starting and terminal vertex of an arrow $a$, respectively. If $g_{ij} = 1$, we also
write $a_{ij}$ instead of $a^{(1)}_{ij} $. 

The quiver $Q$ is called a {\it quiver of  type $C$} and we say the generalized Cartan matrix $C$ is {\it connected} if $Q$ is connected.
Denote by $Q^\circ$ the quiver obtained from $Q$ by deleting all the loops of $Q$. By the condition  $(A2)$, we know that $Q^\circ$ is an acyclic quiver. 

Given an orientation  $\Omega$ of $C$, 
the opposite orientation of  $\Omega$
 is defined as
$\Omega^{op} :=\{(j, i) | (i, j) \in\Omega
\}$.  Denote $\overline{\Omega}=\Omega\cup\Omega^{op}$.
For $(i, j)\in \overline{\Omega}$, 
set
\[\sgn (i,j)=\begin{cases}+1& if\ (i, j) \in \Omega; \\-1&if\ (i, j) \in \Omega^{op}.\end{cases}\]
Let $Q$ be the quiver defined by  $(C,\Omega)$. Let $\overline{Q} = \overline{Q}(C,\Omega)$ be the quiver obtained from $Q$ by adding a new arrow
$a^{(g)}_{ji} : i \to j$ for each arrow $a^{(g)}_{ij} : j \to i \ of\ Q^{\circ}.$

Let $\Omega(i,-)$ be the set $\{j|(i,j)\in\Omega\}.$ Similar one can define ${\Omega}(-,i)$, $\overline{\Omega}(-,i)$ and $\overline{\Omega}(i,-)$.



\begin{definition}
The {\it generalized  preprojective algebra} $\Pi:= \Pi(C,D,\overline{\Omega})$  is the quotient algebra $K\overline{Q}/\overline{I}$ of the path algebra $K\overline{Q}$ by the ideal
 $\overline{I}$ generated by the following relations:

$(P1)$ For each $i\in Q_0$, we have the nilpotency relation
$\varepsilon_i^{c_i} = 0.$

$(P2)$ For each $(i, j) \in\overline{\Omega}$
 and each $1 \leq g \leq g_{ij}$, we have the commutativity relation
$$ a^{(g)}_{ij}\varepsilon_i^{f_{ij}}=\varepsilon_j^{f_{ji}}a^{(g)}_{ij}.$$

$(P3)$ For each $i$, we have the mesh relation
\[\sum_{j\in\overline{\Omega}(-,i)}\sum_{g=1}^{g_{ij}}\sum_{f=0}^{f_{ij}-1}\sgn (i,j)\varepsilon_i^{f_{ij}-1-f}a^{(g)}_{ji}a^{(g)}_{ij} \varepsilon_i^{f}=0.\]
\end{definition}

When the generalized Cartan matrix $C$ is symmetric and the symmetrizer $D$ is the identity matrix, this definition reduces to the one of the classical preprojective algebras for acyclic quivers. 

By the definition of  $\Pi$, it is clear that $\Pi$ does not depend on the orientation of $C$. Hence in the following, we simply write $\Pi=\Pi(C,D)$. 

\begin{rk}\label{r: Pi=C^t}
Our definition of the generalized preprojective algebra here  follows \cite[Definition 7.13]{AHIKM}. This definition is slightly different from the original one given in \cite{GLS17} (which we denote by $\Pi^{\mathrm{GLS}}(C, D)$). As remarked in \cite[Remark 7.14]{AHIKM}, these two definitions describe essentially the same objects. More precisely, the algebra $\Pi(C, D)$ defined in \cite{AHIKM} is isomorphic to the algebra $\Pi^{\mathrm{GLS}}(C^t, D)$ defined in \cite{GLS17}.
\end{rk}

\subsection{The two-sided  ideal $I_i$}

Let $\Pi=\Pi(C,D)$ be a preprojective algebra. For each $i\in Q_0,$ denote by $I_i$ the two-sided ideal $\Pi(1-e_i)\Pi$,  it is easy to see that $e_jI_i=e_j\Pi(1-e_i)\Pi=e_j\Pi$ for $j\neq i$.
Thus we obtain the following decomposition of $I_i$ as right $\Pi$-module 
\begin{eqnarray}\label{g:decomposition of I_i}
    I_i= \bigoplus_{j\in Q_0}e_jI_i= e_iI_i\oplus (\bigoplus_{j\neq i}e_j\Pi).
\end{eqnarray}
Then we have
\begin{lem}\cite{FG19}\label{l:projective resolution of e_iI_i}
  Let $C$ be a symmetrizable Cartan matrix with a symmetrizer $D$ and $\Pi=\Pi(C,D)$. For any $l\leq i\leq n$, 
  \begin{itemize}
      \item[(a)] if  $C$ has no component of Dynkin type, then  $I_i$ is a cofinite tilting ideal of $\Pi$, $\End_{\Pi} I_i\cong\Pi$ and 
      we have the following minimal projective
resolution of $e_iI_i$: \begin{eqnarray}0\ra e_i\Pi\to \bigoplus_{j\in \overline{\Omega}(i,-)}(e_j\Pi)^{-c_{ij}}\to e_iI_i\to 0;\notag\end{eqnarray}.
\item[(b)] if $C$ is of Dynkin type, then $I_i$ is a $\tau$-tilting $\Pi$-module and
we have the following minimal projective
presentation of $e_iI_i$: \begin{eqnarray}e_i\Pi\to \bigoplus_{j\in \overline{\Omega}(i,-)}(e_j\Pi)^{-c_{ij}}\to e_iI_i\to 0.\notag\end{eqnarray}.
  \end{itemize}
  \end{lem}

\begin{theorem}\cite{BIRS09,FG19}\label{t:Non-Dynkin type}
Assume that  $C$ has no component of Dynkin type,  then we have  bijections between the following sets:
\begin{enumerate}
\item[$(1)$] the set of all  cofinite tilting $\Pi$-ideals;
\item[$(2)$] the ideal semigroup $ \langle I_1,I_2,\cdots,I_n\rangle$;
\item[$(3)$] the Weyl group $W(C)$.
\end{enumerate}
Moreover, each $T\in \langle I_1,I_2,\cdots,I_n\rangle$ is a cofinite tilting ideal and $\End_{\Pi}T\cong \Pi$.
\end{theorem}

\begin{theorem}\cite{M14,FG19}\label{t:Dynkin type}
Let $C$ be a symmetrizable Cartan matrix of Dynkin type. Then we have  bijections between the following sets:
\begin{itemize}
\item[$(1)$] the set of all basic support $\tau$-tilting $\Pi$-modules;
\item[$(2)$] the ideal semigroup $ \langle I_1,I_2,\cdots,I_n\rangle$;
\item[$(3)$] the Weyl group $W(C)$.
\end{itemize}
\end{theorem}

\subsection{Geometric representations and  $G$-matrices}\label{s:reflection via g-matrix}
 By \cite{GLS17},  $\Pi=\Pi(C,D)$ is a finite dimensional algebra if and only if $C$ is of Dynkin type. When  $C$ has no component of Dynkin type, for a partial-tilting $\Pi$-module $X$,  we also can define its $g$-vector $g^X$ as subsection~\ref{ss:g-vector}. Hence, for a tilting $\Pi$-module $T$, we can define its $G$-matrix $G_T$.

Let  $w = s_{i_1}s_{i_2}\cdots s_{i_k}$ be a reduced word,  denote by $I_w = I_{i_1}I_{i_2}\cdots I_{i_k}$. If $C$  has no components of Dynkin type, then $I_w$ is a cofinite-tilting ideal. In this case, we fix $G_{I_w}=(g^{e_1I_w},\cdots,g^{e_nI_w})$. If $C$ is of Dynkin type, then $I_w$ is a support $\tau$-tilting $\Pi$ module. Let $(I_w,P_w)$ be the corresponding basic $\tau$-tilting pair, then $P_w=\bigoplus_{1\leq i\leq n,e_iI_w=0}e_{\sigma(i)}\Pi$ (cf. \cite{M14}), where $\sigma \colon Q_0 \to Q_0$ is the Nakayama permutation of $\Pi$.
Let
\begin{align}\label{g:giw}
   g_i(w)=
\begin{cases}
g^{e_i I_w}, & \text{if } e_i I_w \neq 0,\\[4pt]
-\mathbf e_{\sigma(i)} & \text{if } e_i I_w = 0,
\end{cases} 
\end{align}
Then we fix $G_{I_w}=(g_1(w),\cdots,g_n(w))$.  By (\ref{g:decomposition of I_i}) and Lemma~\ref{l:projective resolution of e_iI_i}, it is easy to get that 
\begin{lem}
     Let $C$ be a symmetrizable Cartan matrix with a symmetrizer $D$ and $\Pi=\Pi(C,D)$. Then $G_{I_i}^{t}=R_i$. 
\end{lem}

Moreover, we have  the following general results which is a generalization of \cite[Theorem 6.6]{IR08} for classical preprojective algebras of extended Dynkin type and \cite[Proposition 3.6]{M14} for classical preprojective algebras of Dynkin type
\begin{thm}\label{t: Weyl group via g-matrices}
     Let $C$ be a symmetrizable Cartan matrix with a symmetrizer $D$ and $\Pi=\Pi(C,D)$. Then for any  word $w\in W(C)$, we have $$G_{I_w}=\Sigma_{w^{-1}}=R_w^t.$$
 In particular, the group $\mathcal{G}_{\Pi}$ generated by $\{G_{I_w}^t \mid w\in W(C)\}$ is isomorphic to the Weyl group $W(C)$.    
\end{thm}
    \begin{proof}
    Let us assume first that $C$ has no components of Dynkin type.
Denote by $\varepsilon =\mathcal{K}^b(\operatorname{proj} \Pi)$, {\it i.e.} the homotopy category of bounded complexes of projective $\Pi$-modules. For any $i\in Q_0$, since $I_i$ is a $\Pi$-tilting ideal with $\End_{\Pi}I_i\cong\Pi$, we have an autoequivalence $-\lten_{\Pi}I_i$
of $\varepsilon$  which induces an automorphism $[-\lten_{\Pi}I_i]$
of the Grothendieck group $K_0(\varepsilon)$. Note that $\{[e_j\Pi]~|~j\in Q_0\}$ is a $\mathbb{Z}$-basis of $K_0(\varepsilon)$.
By Lemma~\ref{l:projective resolution of e_iI_i}, 
it is easy to get that for any $1\leq l\leq n$,
$$[e_l\Pi\lten_{\Pi}I_i]=\begin{cases}[e_iI_i]=[\bigoplus\limits_{j\in \overline{\Omega}(i,-)}(e_j\Pi)^{-c_{ij}}]-[e_i\Pi]=[e_i\Pi]-\sum\limits_{j=1}^nc_{ij}[e_j\Pi]& if\ l=i;  \\ [e_lI_i] =[e_l\Pi] & if\ l\neq i, \end{cases}$$
Thus $$[-\lten_{\Pi}I_i]([e_l\Pi])=[e_l\Pi\lten_{\Pi}I_i]=\begin{cases}[e_i\Pi]-\sum\limits_{j=1}^nc_{ij}[e_j\Pi]& if\ l=i;  \\ [e_l\Pi] & if\ l\neq i. \end{cases}$$
Hence, one can get that $$[-\lten_{\Pi}I_i]([e_1\Pi],\cdots,[e_n\Pi])=([e_1\Pi],\cdots,[e_n\Pi])\Sigma_i.$$
Note that for any reduced expression  $w = s_{i_1}s_{i_2}\cdots s_{i_k}$, $I_w = I_{i_1}I_{i_2}\cdots I_{i_k}=I_{i_1}\lten_{\Pi}I_{i_2}\lten_{\Pi}\cdots \lten_{\Pi} I_{i_k}$(cf. \cite{BIRS09,FG19}).
Hence, $[-\lten_{\Pi}I_w]=[-\lten I_{i_1}\lten_{\Pi}I_{i_2}\lten_{\Pi}\cdots \lten_{\Pi} I_{i_k}]=[-\lten_{\Pi}I_{i_k}]\cdots [-\lten_{\Pi}I_{i_1}]$. Thus,
\begin{eqnarray*}
    [-\lten_{\Pi}I_w]([e_1\Pi],\cdots,[e_n\Pi])
    &=&[-\lten_{\Pi}I_{i_k}]\cdots [-\lten_{\Pi}I_{i_1}]([e_1\Pi],\cdots,[e_n\Pi])\\
    &=&([e_1\Pi],\cdots,[e_n\Pi])\Sigma_{i_k}\cdots\Sigma_{i_1}\\
     &=&([e_1\Pi],\cdots,[e_n\Pi])\Sigma_{s_{i_k}\cdots s_{i_1}}\\
    &=&([e_1\Pi],\cdots,[e_n\Pi])\Sigma_{w^{-1}}\\
     &=&([e_1\Pi],\cdots,[e_n\Pi])R_w^t.
\end{eqnarray*}
On the other hand, by $[e_l\Pi\lten_{\Pi}I_w]=[e_lI_w]=([e_1\Pi],\cdots,[e_n\Pi])g^{e_lI_w}$,
we have \begin{eqnarray*}
    [-\lten_{\Pi}I_w]([e_1\Pi],\cdots,[e_n\Pi])
    &=&([e_1\Pi],\cdots,[e_n\Pi])(g^{e_1I_w},\cdots,g^{e_nI_w})\\
    &=&([e_1\Pi],\cdots,[e_n\Pi])G_{I_w}.
\end{eqnarray*}
Because $\{[e_j\Pi]~|~j\in Q_0\}$ is a $\mathbb{Z}$-basis of $K_0(\varepsilon)$,  so we have $$G_{I_w}=\Sigma_{w^{-1}}=R_w^t.$$

If $C$ is of Dynkin type, thus $\Pi$ is a self-injective algebra. When $\Pi$ is a  classical preprojective algebras, Mizuno has proved  $G_{I_w}=\Sigma_{w^{-1}}$ (cf. \cite[Proposition 3.6]{M14}). For the generalized preprojective algebra of Dynkin type,
using same proof as \cite[Proposition 3.6]{M14}, one also can get that $G_{I_w}=\Sigma_{w^{-1}}$. Then, by $\Sigma_{w^{-1}}=R_w^t$, we have $G_{I_w}=\Sigma_{w^{-1}}=R_w^t$.

Since the first part of the theorem yields $G_{I_w}^t = (R_w^t)^t = R_w$, the set $\{G_{I_w}^t \mid w \in W(C)\}$ is precisely the set of geometric representation matrices $\{R_w \mid w \in W(C)\}$. Because the geometric representation is faithful (i.e., an injective group homomorphism), these matrices form a group under multiplication that is isomorphic to $W(C)$.
\end{proof}

\begin{rk}
   For a word $w = s_{i_1}s_{i_2}\cdots s_{i_k}$,  Mizuno denote $\sigma_w^*=\sigma_{i_k}^*\cdots \sigma_{i_1}^*$ in \cite[Definition 3.5]{M14}, where $\sigma_{i}^*$ is just $\rho(s_i)$ of this paper. Then one can get that $\sigma_{w}^*=\sigma^*_{i_k}\cdots \sigma^*_{i_1}=\rho(s_{i_k})\cdots \rho(s_{i_1})=\rho(w^{-1})$.
   Thus $\sigma_{w}^*(\alpha_1^*,\cdots,\alpha_n^*)=\rho(w^{-1})(\alpha_1^*,\cdots,\alpha_n^*)=(\alpha_1^*,\cdots,\alpha_n^*)\Sigma_{w^{-1}}$.
   \end{rk}


    

For a tilting module $T$ over $\Pi$, let $0\ra P_1\ra P_0\ra T\ra 0$ be a minimal projective resolution of $T$. Similar proof as \cite[Proposition 2.5]{AIR14}, one can  get that $P_0$ and $P_1$ have no direct summand in common. Then we have
\begin{cor}
   For any word $w\in W(C)$,  $R_{w}$ has column sign-coherence.
\end{cor}


 Suppose that  $C$ is of Dynkin type. Denote by $$P_{\sigma}=({\bf e}_{\sigma(1)},\cdots,{\bf e}_{\sigma(n)}),$$ where $\sigma \colon Q_0 \to Q_0$ is the Nakayama permutation of $\Pi$. 

\begin{cor} Suppose that  $C$ is of Dynkin type. Let $w,w'$ be two words in $W(C)$, $w_0$  the longest element in $W(C)$. Then
\begin{itemize}
\item[(1)] $G_{I_{w_0}}=-P_{\sigma}$;
    \item[(2)]  $G_{I_{w'}}+P_{\sigma}G_{I_w}=0$ if and  only if $ w'=ww_0$. 
\end{itemize}
   \end{cor}
\begin{proof}
     For $(1)$,  because the support $\tau$-tilting modules have a natural structure
of partially ordered set (cf.~\cite{AIR14}) and $w_0$  the longest element in $W(C)$, then using Theorem~\ref{t:Dynkin type}, one can get that  $I_{w_0}=0$. Thus for each $i$,  $e_iI_{w_0}=0$. By (\ref{g:giw}), we have $g_i(w_0)=-{\bf e}_{\sigma(i)}$. Therefore $$G_{I_{w_0}}=-({\bf e}_{\sigma(1)},\cdots,{\bf e}_{\sigma(n)})=-P_{\sigma}.$$

For $(2)$, using Theorem~\ref{t: Weyl group via g-matrices}, we have
\begin{align*}
w'=ww_0&\iff I_{w'}=I_{ww_0}\\
& \iff  G_{I_{w'}}=G_{I_{ww_0}}=G_{I_{w_0}}G_{I_{w}}=-P_{\sigma}G_{I_w}
   \\
    &\iff G_{I_{w'}}+P_{\sigma}G_{I_w}=0.
\end{align*}
 \end{proof}

\section{Isomorphism induced by a 2-term silting complex  and the G-matrix}\label{s: iso induced by 2-term silting complex and G-matrix}
\subsection{2-term silting complex}
Let $A$ be a finite dimensional algebra. 
We consider 2-term  complexes \( P^* \) in the bounded homotopy category \( K^b(\mathsf{proj}\, A) \). These are complexes \( P^* = \{P^i\} \) with \( P^i = 0 \) for \( i \neq -1, 0 \). Such a complex is called \emph{pre-silting} if \( \mathrm{Hom}_{K^b(\mathsf{proj}\, A)}(P^*, P^*[1]) = 0 \) and \emph{silting} if in addition \( \mathsf{thick}\, P = K^b(\mathsf{proj}\, A) \). 

We denote by $\mathsf{2\text{-}silt}\, A$ (respectively, $\mathsf{2\text{-}presilt}\, A$) the set of isomorphism classes of basic two-term silting (respectively, presilting) complexes for $A$.

\begin{lem}\cite[Theorem 3.2]{AIR14}
     Let \( A \) be a finite dimensional \( k \)-algebra. Then there exists a bijection 
\[
\phi:\mathsf{2\text{-}silt} A \longleftrightarrow \mathsf{s}\tau\mathsf{\text{-}tilt}  A
\]
such that $\phi((P^{1}\xra{d} P^0))= H^0(P^*)$ and $\phi^{-1}((M,P))=(P_1 \oplus P\xra{(f,\ 0)} P_0)$
where \( f: P_1 \to P_0 \) is a minimal projective presentation of \( M \). Moreover, $\phi$ induce a bijection between $\mathsf{2\text{-}presilt}\, A$ and the set of basic $\tau$-rigid $A$-pairs.
\end{lem}

Let $P^*=(P^{1}\xra{d} P^0)$ be a basic 2-term silting complex in $K^b(\proj A)$. 
Suppose that $P^*=\bigoplus_{j=1}^nP^*_j$, where each $P^*_j=(\bigoplus_{i=1}^nP(i)^{b_{ij}}\ra \bigoplus_{i=1}^nP(i)^{a_{ij}})$ is indecomposable. It is easy to get that the matrix $(a_{ij}-b_{ij})_{1\leq i,j\leq n}$ is a $G$-matrix of $\phi(P^*)$, which is also called  a {\it $G$-matrix of $P^*$},  simply denoted by $G_{P^*}$.

Let $B=\End_{K^b(\proj A)} P^*$ (for more results on $\mod B$, we refer to \cite{BZ16}). By \cite{Hu22},  the correspondence 
      $$f:\dimv X\to \dimv \Hom_{\md^b(A)}(P^*,X)-\dimv \Hom_{\md^b(A)}(P^*,X[1]),$$
      where $X$ is an $A$-module, induces an isomorphism $f:K_0(A)\to K_0(B)$ of the Grothendieck groups of $A$ and $B$.   Using Lemma~\ref{l:key lemma}, one can easily obtain the following result.
\begin{thm}\label{t:iso induced by a 2-term silting complex}
For each $A$-module $X$,
$$G_{P^*}^t\dimv X= \dimv \Hom_{\md^b(A)}(P^*,X)-\dimv \Hom_{\md^b(A)}(P^*,X[1]).$$
\end{thm}
\begin{proof} Here, we provide a direct proof.

Note that for each simple $A$-module $S_l$ for  $1\leq l\leq n$,
\begin{eqnarray*}
   f(\dimv S_l)&=&\dimv \Hom_{\md^b(A)}(P^*,S_l)-\dimv \Hom_{\md^b(A)}(P^*,S_l[1])\\
   &=&(\dim \Hom_B(\Hom_{\md^b(A)}(P^*,P_j^*),\Hom_{\md^b(A)}(P^*,S_l))_{1\leq j \leq n }\\
   &&-(\dim \Hom_B(\Hom_{\md^b(A)}(P^*,P_j^*),\Hom_{\md^b(A)}(P^*,S_l[1]))_{1\leq j \leq n }\\
 &=&(\dim\Hom_{\md^b(A)}(P_j^*,S_l))_{1\leq j\leq n}-(\dim\Hom_{\md^b(A)}(P_j^*,S_l[1]))_{1\leq j\leq n}\\
   &=&(\dim \Hom_{A}( \bigoplus_{i=1}^nP(i)^{a_{ij}},S_l))_{1\leq j\leq n}-(\dim \Hom_{A}( \bigoplus_{i=1}^nP(i)^{b_{ij}},S_l))_{1\leq j\leq n}
   \\&=&(a_{lj}-b_{lj})_{1\leq j\leq n},
\end{eqnarray*}
which is the $l$-th column of $G_{P^*}^t$.
Thus $$(f(S_1),\cdots,f(S_n))=G^t_{P^*}.$$
Therefore, for each $A$-module $X$,
$$G_{P^*}^t\dimv X= \dimv \Hom_{\md^b(A)}(P^*,X)-\dimv \Hom_{\md^b(A)}(P^*,X[1]).$$
\end{proof}

\noindent\textbf{Acknowledgements}
The author thanks Changjian Fu and Yu Zhou for helpful discussions.

\noindent\textbf{Consent to Publish} 
Not applicable.

\noindent\textbf{Author Contributions}  I am the only author who contributed to this research.

\vspace{0.5em}
\noindent\textbf{Data Availability} No datasets were generated or analysed during the current study.

\vspace{0.8em}
\noindent{\Large\textbf{Declarations}}

\vspace{0.8em}
\noindent\textbf{Competing interests} The authors declare no competing interests.

\vspace{0.5em}
\noindent\textbf{Ethical Approval} Not applicable.

\end{document}